\input amstex
\documentstyle{amsppt}
\magnification 1200
\NoRunningHeads
\NoBlackBoxes
\document

\def\tr{{\text{tr}}}
\def\ell{{\text{ell}}}
\def\Ad{\text{Ad}}
\def\u{\bold u}
\def\m{\frak m}
\def\O{\Cal O}

\def\qdet{\text{qdet}}

\def\RR{\Bbb R}

\def\h1{\hat{\bold 1}}

\def\deg{\text{deg}}

\def\pF{F^{\prime}}

\def\bg{\overline{\g}}
\def\bG{\overline{G}}

\def\ad{\text{ad}}
\def\Hom{\text{Hom}}
\def\hh{\hat\h}
\def\a{\frak a}

\def\Ua{U_q(\tilde\g)}
\def\U2{{\Ua}_2}
\def\g{\frak g}
\def\n{\frak n}
\def\hh{\frak h}

\def\Z{\Bbb Z}

\def\d{\partial}

\def\l{\lambda}

\def\z{\bold z}
\def\Id{\text{Id}}
\def\<{\langle}
\def\>{\rangle}
\def\o{\otimes}
\def\e{\varepsilon}

\def\Id{\text{Id}}
\def\End{\text{End}}

\def\b{{\frak b}}

\topmatter
\title Quantization of Lie bialgebras, III
\endtitle
\author {\rm {\bf Pavel Etingof and David Kazhdan} \linebreak
\vskip .1in
Department of Mathematics\linebreak
Harvard University\linebreak 
Cambridge, MA 02138, USA\linebreak
e-mail: etingof\@math.harvard.edu\linebreak kazhdan\@math.harvard.edu}
\endauthor
\endtopmatter

\head Introduction\endhead

This article is the third part of the series of papers on quantization of Lie 
bialgebras which we started in 1995. However, its object of study is
much less general than in the previous two parts. While in the first 
and second paper we deal with an arbitrary Lie bialgebra, here we study 
Lie bialgebras of $\g$-valued functions on a punctured rational or elliptic 
curve, where $\g$ is a finite dimensional simple Lie algebra. Of course, the
general result of the first paper, which says that any Lie bialgebra 
admits a quantization, applies to this particular case. However, this result 
is not sufficiently effective, as the construction of quantization 
utilizes a Lie associator, which is computationally unmanageable.
The goal of this paper is to give a more effective quantization procedure
for Lie bialgebras associated to punctured curves, i.e. a procedure
which will not use an associator. We will describe a 
general quantization procedure
which reduces the problem of quantization of the algebra of 
$\g$-valued functions 
on a curve with many punctures to the case of one puncture, and apply
this procedure in a few special cases to obtain an explicit 
quantization. 

The main object of study in this paper are Lie bialgebras associated 
to rational and elliptic curves with punctures, which can be described 
as follows. 

We work over an algebraically closed field $k$ of characteristic zero.
 Let $\Sigma$ be a 1-dimensional algebraic group over $k$
(i.e. $\Bbb G_a,\Bbb G_m$, or an elliptic curve), and $u$ be an additive 
formal parameter near the origin. Let 
$r\in \g\o\g(\Sigma)$ be a rational $\g\o\g$-valued function on $\Sigma$ 
with the Laurent expansion near $0$ of the form $\sum_\alpha X_\alpha
\o X_\alpha^*/u+O(1)$,
satisfying the classical Yang-Baxter equation (2.1)
(here $\{X_\alpha\},\{X_\alpha^*\}$ are dual bases of $\g$ with respect to the 
half-Killing form). 
Such a function is called a classical r-matrix. 

Consider the vector space $\a=t^{-1}\g[t^{-1}]$. 
Let $\tau(u)=\sum_{\alpha}\sum_{m\ge 1} (X_\alpha\o X_\alpha^*
t^{-m})u^{m-1}\in\g\o\a[[u]]$. 
Define a Lie bialgebra structure on $\a$ by the formulas
$$
\gather
[\tau^{13}(u),\tau^{23}(v)]=
[r^{12}(u-v),\tau^{13}(u)+\tau^{23}(v)],\\
\delta(\tau(u))=[\tau^{12}(u),\tau^{13}(u)].\tag 1
\endgather
$$ 

It is convenient to understand the classical r-matrix 
as a bilinear form $\beta$ on $\a$ with values in 
$k((u))$, defined by the rule $\beta(x,y)(u)=-\text{Res}_{v,w=0}
\<x(v)\o y(w),r(v-w+u)\>$, where $\<,\>$ is
the invariant form on $\g$, and $r(v-w+u)$ is understood as an element
of $\g\o\g((u))[[v,w]]$ (that is, the function $\frac{1}{v-w+u}$ is expanded 
as $\sum u^{-r-1}(w-v)^r$). We can regard $\beta$ as an element of 
$\a^*\o\a^*((u))$, where the tensor product is understood in the completed 
sense. The form $\beta$ satisfies a set of axioms which are dual 
to the axioms of a pseudotriangular structure on a Lie bialgebra \cite{Dr1}. 
Thus we call $\beta$ a copseudotriangular structure. 
 
Now we describe Lie bialgebras  corresponding to punctured curves. 
Let $\Gamma$ be the set of poles of $r$.
For any collection of points $\z=(z_1,...,z_n)$ on $\Sigma$
such that $z_i-z_j\notin\Gamma$ when $i\ne j$,
define the Lie bialgebra $\a_\z$ to be the direct sum
$\a_\z=\a_{z_1}\oplus...\oplus \a_{z_n}$, where $\a_{z_i}$ are Lie 
subbialgebras in $\a_\z$ identified with $\a$, and the commutation 
relations between $\a_{z_i}$ and $\a_{z_j}$ are given by the formula
$$
[\tau^{13}_i(u),\tau^{23}_j(v)]=
[r^{12}(u-v+z_i-z_j),\tau^{13}_i(u)+\tau^{23}_j(v)],\tag 2
$$ 
where $\tau_i$ is the image of $\tau$ under the identification
$\a\to\a_{z_i}$. Here $r(u-v+z_i-z_j)$ is regarded as an element of
$\g\o \g[[u,v]]$. This is possible as $r$ is regular at $z_i-z_j$ 
by the choice of $\z$. 

The simplest example of this situation occurs when $\Sigma=\Bbb G_a$
and $r(u)=\sum_i X_\alpha\o X_\alpha^*/u$ (the Yang's r-matrix). 
In this case it is easy to check
that $\a$ is the Lie algebra $t^{-1}\g[t^{-1}]$, with cobracket
dual to the standard bracket in $\g[[t]]$ \cite{Dr1}, and $\a_\z$ is the
Lie algebra
of rational functions on the line with values 
in $\g$ which have no poles outside $z_1,...,z_n$ and vanish at infinity,
with a Lie coalgebra structure dual to the standard bracket
in $\g[[t_1]]\oplus...\oplus\g[[t_n]]$ \cite{Dr1}.

The Lie bialgebra $\a_\z$ has two essential 
properties:

1. {\it Local factorization in a product of equal components}. 

(a) {\it factorization}: As a Lie coalgebra,
$\a_\z$ admits a decomposition $\a_\z=\a_1\oplus...\oplus\a_n$,
where $\a_i$ are Lie subbialgebras of $\a_\z$.

(b) {\it locality}: $[\a_i,\a_j]\subset 
\a_i\oplus\a_j$.

(c) {\it equal components}: All the Lie bialgebras $\a_i$ are identified
with the same Lie bialgebra $\a$. 

2.  {\it Expression of commutator via the classical r-matrix}.   

The commutator between the components $\a_i$ and $\a_j$ is given in 
terms of the classical r-matrix via formula (2). 

Thus, we see that the structure of $\a_\z$ is completely determined by 
the pair $(\a,\beta)$ of a Lie bialgebra and a copseudotriangular structure 
on it. 

Our purpose in this paper is to describe an ``explicit'' 
quantization of $\a_\z$. The construction of such quantization consists 
of two parts. 

Part 1. We define the notion of a copseudotriangular structure $B$ 
on a Hopf algebra $A$, by dualizing the notion of a pseudotriangular 
structure, introduced by Drinfeld \cite{Dr1}. 
We show that any nondegenerate copseudotriangular 
Lie bialgebra $(\a,\beta)$ can be quantized, i.e. that there exists 
a copseudotriangular Hopf algebra $(A,B)$ whose quasiclassical limit is 
$(\a,\beta)$. This is done using the methods of \cite{EK1,EK2}
 
Part 2. We define the notion of a factored Hopf algebra, 
which is a quantization of the above notion of a 
factored Lie bialgebra, and give a construction of a factored Hopf algebra 
$A_\z$ corresponding to points $(z_1,...,z_n)$ starting from a 
copseudotriangular Hopf algebra $(A,B)$. We show that if $(A,B)$ is a 
quantization of $(\a,\beta)$ then $A_\z$ is a quantization of $\a_\z$.  

In some cases, this quantization of $\a_\z$ can be described 
by explicit formulas. For example, 
in the case of the Yang's r-matrix the quantization of $(\a,\beta)$ 
is $(Y^*(\g),\Bbb R^{-1})$, where $Y^*(\g)$ is the dual algebra to the Yangian 
$Y(\g)$, with opposite product, 
and $\Bbb R\in Y(\g)\o Y(\g)((u))$ is the pseudotriangular
structure on $Y(\g)$ defined by Drinfeld. 
In this case, $A_\z$ is defined by (3.5),(3.6), where 
$R$ is the Yang's quantum R-matrix. 

Let us describe briefly the structure of the paper. 

In Chapter I we introduce the notion of a  
locally factored Lie bialgebra, which is a Lie bialgebra with
properties 1(a) and 1(b), and the corresponding 
quantum notion of a locally factored Hopf algebra.
We also introduce the notion of a weak
($\d$)-copseudotriangular Lie bialgebra,
which is a Lie bialgebra $\a$ with a derivation $\d$ 
and a form $\beta:\a\o\a\to k((u))$ 
having analogous properties to the form $\beta$ above, and the corresponding
quantum notion of a weakly 
($\d$)-copseudotriangular Hopf algebra. We explain how 
to introduce a factored Lie bialgebra structure on $\a^{\oplus n}$ 
given a weak $\d$-copseudotriangular 
structure $\beta$ on a Lie bialgebra $\a$, and how 
to introduce a factored Hopf algebra structure on $A^{\o n}$ 
given a weak $\d$-copseudotriangular structure $B$ on a Hopf algebra $A$.
Then we explain how to quantize a given weak $\d$-copseudotriangular 
structure on a Lie bialgebra. This allows us to describe the quantization
(constructed in \cite{EK1,EK2}) 
of the factored algebra $\a^{\oplus n}$ purely in terms of the quantization 
of $(\a,\beta)$. 

In Chapter II we describe in detail the Lie bialgebras $\a$ and $\a_\z$ 
discussed above, and consider various examples (rational, trigonometric,
and elliptic). The results of this Chapter are not original.

In Chapter III we give the main results of the paper. 
In the first part we completely describe the quantization 
$U_{A_\z}$ of $\a_\z$ in the
case of the Yang's r-matrix in terms of the Yangian $Y(\g)$ and its 
universal R-matrix. In the second part we consider the case 
$\g=\frak{sl}_N$ and $\frak{gl}_N$, and describe quantizations of the
simplest rational, trigonometric, and elliptic algebras $\a_\z$ 
explicitly by generators and relations containing 
the quantum R-matrices $R(u)\in\End(k^N\o k^N)[[h]]$. 
This construction is analogous to the Faddeev-Reshetikhin-Takhtajan
construction of quantum $GL_N$ \cite{FRT}, and 
the similar construction of the Yangian of $\frak{gl}_N$ \cite{Dr1}. 

{\bf Remark.} Algebras similar to $U_{A_\z}$ were considered 
in \cite{Ch1} and \cite{AGS1,AGS2}.

In a subsequent paper we will use the algebras $U_{A_\z}$ to
describe a deformed version of
the holomorphic part of the Wess-Zumino-Witten conformal field theory
(in genus zero). More specifically, we will
consider representation theory of
$U_{A_\z}$, the notion of quantum conformal blocks, 
quantum vertex operators, and show how the quantum Knizhnik-Zamolodchikov
equations of Frenkel-Reshetikhin appear naturally in this context.

\head Acknowledgements\endhead

The authors would like to thank I.Cherednik for useful discussions 
and references, and V.Schomerus for careful reading of the manuscript. 
The authors are grateful to V.Drinfeld who pointed out 
an error in the original formulation of Proposition 3.9. 
The work of P.E. was supported by an NSF postdoctoral 
fellowship. The work of D.K. was supported by an NSF grant.

\head 1. Copseudotriangular Lie bialgebras
and their quantization\endhead

Throughout the paper, $k$ will denote an algebraically closed
field of characteristic zero, and ``an algebra'' means ``an associative 
$k$-algebra with 1''.

\subhead 1.1. Factored algebras.\endsubhead

Let $A$ be an algebra over $k$.

\proclaim{Definition} A factorization of $A$ is a collection
of algebras $A_1,...,A_n$, which are subalgebras
in $A$, such that for any 
$\sigma\in S_n$ the multiplication map defines a bijection
$A_{\sigma(1)}\o...\o A_{\sigma(n)}\to A$. 
Such a factorization is called local if the image of
$A_i\o A_j$ in $A$ under the multiplication map
is a subalgebra of $A$ for all $i,j$.
\endproclaim

It is clear that a factorization of $A$ is local if and only if 
$A_iA_j=A_jA_i$ for all $i,j$. 

We call an algebra $A$ equipped with a (local) factorization a (locally) 
factored algebra. 

Let $A_1,...,A_n$ be algebras.
We want to describe 
locally factored algebras with factors $A_1,...,A_n$.

Suppose we are given a locally factored algebra $A$,
with factors $A_1,...,A_n$, $n\ge 2$, 
and product 
$m:A\o A\to A$. Let $m_{ij}: A_i\o A_j\to A$ be the restriction of $m$ to 
$A_i\o A_j$. By the definition, the map $m_{ij}$ is 
injective for $i\ne j$, and 
the images of $m_{ij}$ and $m_{ji}$ are the same. Therefore, 
for any $i<j$, $i,j\in\{1,...,n\}$, we can define a linear isomorphism
$X_{ij}: A_i\o A_j\to A_i\o A_j$ by the formula 
$X_{ij}=m_{ij}^{-1}\circ m_{ji}\circ \sigma$, where
$\sigma$ is the permutation of components.

It is easy to see that the maps $X_{ij}$ satisfy the following 
conditions:

(i) 
$$
X_{ij}X_{ik}X_{jk}=X_{jk}X_{ik}X_{ij}, i<j<k;\tag 1.1
$$

(ii) 
$$
X_{ij}(a\o bc)=(1\o m)(X^{12}_{ij}X^{13}_{ij}(a\o b\o c)),\
X_{ij}(ab\o c)=(m\o 1)(X^{23}_{ij}X^{13}_{ij}(a\o b\o c)),\tag 1.2
$$

(iii)
$$
 X_{ij}(a\o 1)=a\o 1, a\in A_i,\ X_{ij}(1\o a)=1\o a, a\in A_j.\tag 1.3
$$

{\bf Remark.} Condition (i) is vacuous if $n=2$.
 
Conversely, let $X_{ij}: A_i\o A_j\to A_i\o A_j$ be an arbitrary 
collection of invertible operators satisfying (i)-(iii), and
$I(\{X_{ij}\})$ be the ideal in the free product
$A_1*...*A_n$ 
generated by the relations
$ba=m(X_{ij}(a\o b))$, $a\in A_i,b\in A_j$, $i<j$.
Denote by $P_n(A_1,...,A_n,\{X_{ij}\})$ 
the quotient of $A_1*...*A_n$ by this ideal. 

Let $\phi: A_1\o...\o A_n\to P_n(A_1,...,A_n,\{X_{ij}\})$ 
be the linear map defined by
the rule $a_1\o...\o a_n\to a_1...a_n$. 

\proclaim{Proposition 1.1} The map $\phi$ is a bijection. 
\endproclaim

To prove this proposition, we
define a product on the space $A_1\o...\o A_n$ by 
$$
\gather
(a_1\o...\o a_n)\cdot (b_1\o...\o b_n)=\\
(m\o...\o m)(X_{n-1,n}^{2n-1,n}...X_{23}^{n+2,3}...X_{2n}^{n+2,n}X_{12}^{n+1,2}...X_{1n}^{n+1,n}(a_1\o...\o a_n\o b_1\o...\o b_n)),\tag 1.4
\endgather
$$
where the superscripts denote the components where the corresponding 
$X$-operator acts.

\proclaim{Lemma 1.2} The product defined by (1.4)
endows $A_1\o...\o A_n$ with the structure of an associative algebra,
in which the unit is $1\o...\o 1$. 
\endproclaim

\demo{Proof} The associativity follows directly from properties (i),(ii)
of $X_{ij}$. The unit axiom follows from property (iii) of $X_{ij}$.
$\square$\enddemo

We denote the algebra of Lemma 1.2 by $E$. 
We have a natural
homomorphism $\psi: P_n(A_1,...,A_n,\{X_{ij}\})\to E$, 
induced by the embeddings $A_i\to E$. 

\proclaim{Lemma 1.3} $\psi$ is an isomorphism. 
\endproclaim

\demo{Proof}
This homomorphism is surjective, as 
$E$ is generated by $A_i$. It is also injective,
as $P_n(A_1,...,A_n,\{X_{ij}\})$ is obviously 
spanned by the elements of the form 
$a_1...a_n$, 
$a_i\in A_i$. The Lemma is proved.
$\square$\enddemo
 
\demo{Prooof of Proposition 1.1} It is clear that $\psi\circ \phi=\Id$, so
$\phi$ is an isomorphism. 
$\square$\enddemo

\proclaim{Corollary 1.4} The assignment $A\to \{X_{ij}(A)\}$ is
a 1-1 correspondence between locally factored algebra structures
on $A_1\o...\o A_n$,  with factors
$A_1,...,A_n$, and collections of
invertible operators $\{X_{ij}\}$ satisfying (i)-(iii).
\endproclaim

The same definitions and results apply to the case when $A$ is an algebra 
over $k[[h]]$ which is a deformation of an algebra over $k$. In this case, 
the sign $\o$ denotes the completed (in the h-adic topology)
tensor product over $k[[h]]$.  
 
\subhead 1.2. Factored Lie bialgebras and Hopf algebras\endsubhead

Let $\a$ be a Lie bialgebra over $k$.

\proclaim{Definition} A factorization of $\a$ is a decomposition
$\a=\a_1\oplus...\oplus\a_n$, where $\a_i\subset \a$ are Lie subbialgebras. 
Such a factorization is called local if 
$\a_i\oplus \a_j$ is a Lie subbialgebra of $\a$ for $i\ne j$, i.e. 
if $[\a_i,\a_j]\subset \a_i\oplus\a_j$ 
for $i\ne j$.
\endproclaim

Let $U$ be a Hopf algebra.

\proclaim{Definition} A factorization of $U$ is a collection
of Hopf algebras $U_1,...,U_n$, which are Hopf subalgebras
in $U$, such that for any $\sigma\in S_n$ 
the multiplication map defines a bijection
$U_{\sigma(1)}\o...\o U_{\sigma(n)}\to U$. 
Such a factorization is called local if the image of
$U_i\o U_j$ in $U$ under the multiplication map
is a Hopf subalgebra of $U$ for $i\ne j$, i.e. 
if $U_jU_i=U_iU_j$ 
for $i\ne j$.
\endproclaim

We will call a Lie bialgebra (Hopf algebra) with a (local) factorization
a (locally) factored Lie bialgebra (Hopf algebra).

It is clear that if a Lie bialgebra $\a$ is the quasiclassical limit of
a quantized universal enveloping (QUE) algebra $U$ 
(\cite{Dr1}), then any factorization of $U$ 
(as a QUE algebra) defines a factorization
of $\a$. Moreover, if the first factorization is local, so is the second one. 

Fix a universal Lie associator $\Phi$.
Let $\a\to U_h(\a)$ be the functor of quantization from the category
of Lie bialgebras to the category of QUE algebras, defined in \cite{EK2}
using $\Phi$.
Suppose we are given a factorization $\a=\oplus\a_i$.
Using the functoriality of $U_h$, we obtain embeddings of Hopf algebras
$U_h(\a_i)\to U_h(\a)$, so we can regard $U_h(\a_i)$ as Hopf subalgebras
of $U_h(\a)$. These Hopf subalgebras define a factorization
$U_h(\a)=U_h(\a_1)\o...\o U_h(\a_n)$.
Moreover, it follows from the functoriality
of $U_h$ that if the factorization of $\a$ is local, so is the factorization
of $U_h(\a)$.

Thus, we have the following proposition.

\proclaim{Proposition 1.5} To any (local) factorization of
a Lie bialgebra $\a$ one can assign its quantization, i.e.
a (local) factorization of $U_h(\a)$, whose
quasicalssical limit is the initial factorization of $\a$.
\endproclaim

Now we give an analogue of Corollary 1.4 for Hopf algebras.

Let $U_1,...,U_n$ be Hopf algebras. Then the free product $U_1*...*U_n$
has a natural Hopf algebra structure induced by the Hopf algebra 
structures on $U_1,...,U_n$.  

\proclaim{Proposition 1.6} The assignment 
$U\to\{X_{ij}(U)\}$ is a 1-1 correspondence
between locally factored Hopf algebra structures on $U_1\o...\o U_n$,
 with factors $U_1,...,U_n$, 
and such collections of
invertible operators $\{X_{ij}\}$ satisfying (i)-(iii), that 
the ideal $I(\{X_{ij}\})$ is a Hopf ideal in $U_1*...*U_n$. 
\endproclaim

\demo{Proof} Clear.
$\square$\enddemo

\subhead 1.3. Copseudotriangular Lie bialgebras\endsubhead

Let $\b$ be a Lie bialgebra over $k$, 
with commutator $\mu$ and cocommutator $\delta$. 

Recall \cite{Dr1} that a quasitriangular 
structure on $\b$ is an element $r\in \b\o \b$ satisfying 
the classical Yang-Baxter equation 
$$
[r^{12},r^{13}]+[r^{12},r^{23}]+[r^{13},r^{23}]=0,\tag 1.5
$$
such that 
$$
\delta(x)=[x\o 1+1\o x,r], x\in\b.\tag 1.6
$$ 

It is easy to see that relation (1.5) can be 
replaced by any one of the two relations
$$
(\delta\o 1)(r)=[r^{13},r^{23}], 
(1\o\delta)(r)=[r^{13},r^{12}].\tag 1.7 
$$
Indeed, given (1.6), relation (1.5) is equivalent to either of (1.7).

Let $\b$ be a finite dimensional quasitriangular Lie bialgebra, and $\a=\b^*$.
In this case, the element $r$ 
defines a bilinear form 
$\beta: \a\o\a\to k$, by
$$
\beta(x,y)=(r,x\o y),\ x,y\in\a.\tag 1.8
$$
This form has the following properties, which are equivalent to
properties (1.6),(1.7) of $r$: 

$$
\beta([xy],z)=\beta(x\o y,\delta(z)), \beta(x,[zy])=\beta(\delta(x),y\o z);
\tag 1.9
$$

$$
[yx]=\beta_{12}(x\o\delta(y))+\beta_{13}
(\delta(x)\o y).
\tag 1.10
$$

\proclaim{Definition} A bilinear form $\beta$ on a Lie bialgebra $\a$ 
satisfying (1.9)-(1.10) is called a coquasitriangular structure on $\a$. 
\endproclaim

Thus, if $\a$ is finite-dimensional, 
a coquasitriangular structure on $\a$ is the same thing as a quasitriangular
structure on $\a^*$. 

The following generalization of this definition is 
essentially due to Drinfeld, \cite{Dr1}.

Let $\a$ be a Lie bialgebra over $k$.
Let $\d$ be a derivation of $\a$ as a Lie bialgebra.
Define the map $\alpha_u:\a\to\a[[u]]$ by 
$$
\alpha_u(x)=e^{u\d}x.\tag 1.11
$$ 

\proclaim{Definition} 
A $\d$-copseudotriangular structure on
a Lie bialgebra $\a$ is a bilinear form
 $$
\beta:\a\o\a\to k((u)),\tag 1.12 
$$
which satisfies the following
 conditions:

$$
\beta([xy],z)=\beta(x\o y,\delta(z)), \beta(x,[zy])=\beta(\delta(x),y\o z).
\tag 1.13
$$

$$
[y,\alpha_{u}(x)]=\beta_{12}(x\o\delta(y))(u)+\alpha_{u}(\beta_{13}
(\delta(x)\o y))(u).
\tag 1.14
$$

$$
\beta(\alpha_v(x),y)(u)=\beta(x,\alpha_{-v}(y))(u)=\beta(x,y)(u+v).\tag 1.15 
$$

If $\O\subset k((u))$ is a subalgebra, and in addition to (1.13)-(1.15)
$\beta$ takes values in $\O$, we say
that $\beta$ is a $\d$-copseudotriangular $\O$-structure. 
\endproclaim

Consider the linear map
$d=-\frac{1}{2}\mu\circ \delta:\a\to \a$. It is easy
to check that this map is a derivation of $\a$ as a Lie bialgebra \cite{Dr1}. 
We call it the canonical derivation
of $\a$. 
If $\d=d$, we will call $\beta$ a copseudotriangular structure
(without specification of $\d$).

{\bf Example.} If $\O=k$ then a $\d$-copseudotriangular $\O$-structure is 
the same thing as a coquasitriangular structure, which vanishes on the image
of $\d$.

\proclaim{Proposition 1.7} Let $\beta$ be a
$\d$-copseudotriangular structure on $\a$ 
such that for a suitable $m\in\Z$ 
we have $\beta(x,y)\in u^mk[[u]]$ for any $x,y\in\a$. Then 
$\beta$ is constant (does not depend on $u$), 
and therefore defines a coquasitriangular structure on $\a$, which vanishes
on the image of $\d$. 
\endproclaim

\demo{Proof} Consider the form
$\beta_0:\a\o\a\to k((u))$ defined by
$\beta_0(x,y)(u):=\beta(\alpha_{-u}(x),y)(u)$. Since
$\beta$ takes values in $u^mk[[u]]$, this expression makes sense. 
From (1.15) we get that $\beta_0(u+v)=\beta_0(u)$, so $\beta_0$ is constant.
Identities (1.13)-(1.14) for $\beta$ imply (1.9)-(1.10) for 
$\beta_0$. Thus, 
$\beta_0$ is a quasitriangular structure on $\a$.
It is obvious that $\beta_0$ vanishes on the image of $\d$.
$\square$\enddemo

{\bf Remark.} Proposition 1.7 shows that interesting (i.e. not coquasitriangular)
examples of $\d$-copseudotriangular structures can only arise  
when $\beta(x,y)(u)$ can have a pole of arbitrary order
at $u=0$, for which $\a$ has to be infinite dimensional.

\subhead 1.4. Copseudotriangular Hopf algebras\endsubhead

Let $A$ be a Hopf algebra over $k[[h]]$, such that 
$A_0=A/hA$ is commutative, and 
$A$ is a deformation of $A_0$. Let $m,1,\Delta,\e,S$ be the product, unit,
coproduct, counit, and the antipode of $A$. 

Let $\d$ be a derivation of $A$. 
 Define a Hopf algebra homomorphism
$$
\alpha_u:=e^{u\d}: A\to A[[u]].\tag 1.16
$$

\proclaim{Definition} A $k[[h]]$-bilinear 
form $B:A\o A\to k((u))[[h]]$ is called a $\d$-copseudotriangular structure
on $A$ if it satisfies the following conditions:
$$
B(xy,z)=B(x\o y,\Delta(z)), 
B(x,zy)=B(\Delta(x),y\o z), x,y,z\in A;\tag 1.17
$$

$$
m((\alpha_{u}\o 1)[B_{13}(\Delta(x)\o\Delta(y))(u)])=
m^{op}((\alpha_{u}\o 1)[B_{24}(\Delta(x)\o\Delta(y))(u)]);\tag 1.18
$$

$$
B(\alpha_v(x),y)(u)=B(x,\alpha_{-v}(y))(u)=B(x,y)(u+v);\tag 1.19
$$

$$
B(1,x)=B(x,1)=\e(x), x\in A; B(x,y)=\e(x)\e(y)+O(h),\tag 1.20
$$ 
where $B_{ij}$ means that $B$ is evaluated on the $i$-th and $j$-th 
components of the product.

Let $\O$ be a subalgebra of $k((u))$. We say that $B$ is
a $\d$-copseudotriangular
$\O$-structure on $A$ if $B$ takes values in
$\O[[h]]$. 
\endproclaim

{\bf Remark.} Commutativity of $A_0$ is essential to enable 
the equality $B(x,y)=\e(x)\e(y)+O(h)$ in presence of condition 
(1.18).

Since $A_0$ is commutative, 
we have $S^2=1+O(h)$. Let
$D:=\frac{1}{h}\ln S^2$. It is clear that $D$ is a derivation of $A$.
We call $D$ the canonical derivation
of $A$. If $\d=D$, we will call $B$ a copseudotriangular structure
(without specification of $\d$).

\subhead 1.5. Weak $\d$-copseudotriangular structures\endsubhead

We will need the following
 weaker version of the notion of a $\d$-copseudotriangular 
structure.

\proclaim{Definition} (i) Let $\a$ be a Lie bialgebra over $k$.
A $k((u))$-valued bilinear form $\beta$ on $\a$ 
is called a weak $\d$-copseudotriangular structure on 
$\a$ if it satisfies equations (1.13),(1.15), and 
the equation
$$
\beta([y,\alpha_{u}(x)],z)(v)=\beta_{12}(u)\beta_{34}(v)(x\o\delta(y)\o z)-
\beta_{14}(v)\beta_{23}(u)
((\alpha_{u}\o 1)(\delta(x))\o y\o z).
\tag 1.21
$$

(ii) Let $A$ be a Hopf algebra, as in Section 1.4. 
A $k((u))[[h]]$-valued bilinear form $B$ on $A$ 
is called a weak $\d$-copseudotriangular structure on 
$A$ if it satisfies equations (1.17),(1.19),(1.20) and 
the equation
$$
B\left( m((\alpha_{u}\o 1)[B_{13}(\Delta(x)\o\Delta(y))(u)]),z\right)(v)=
B\left( m^{op}((\alpha_{u}\o 1)[B_{24}(\Delta(x)\o\Delta(y))(u)]),z\right)(v).
\tag 1.22
$$

If $\beta$ or $B$ takes values in $\O\subset k((u))$, it is called a
weak $\d$-copseudotriangular $\O$-structure. 
\endproclaim

If $\d=d$ (respectively, $\d=D$), we will call $\beta$ (respectively, $B$)
 a weak copseudotriangular structure (without specification of $\d$).

{\bf Remark.} Equations (1.21), (1.22) are obtained 
by applying the functionals $\beta(*,z)$, $B(*,z)$ to equations (1.14),(1.18).
Thus,  for a left-nondegenerate bilinear form
(i.e. a form with trivial left kernel), the property 
to be a weak $\d$-copseudotriangular structure
is the same as to be a $\d$-copseudotriangular structure.

\subhead 1.6.  $\d$-copseudotriangular structure on h-formal groups\endsubhead

Let $\a_0$ be a Lie coalgebra, $A_0=\prod_{j\ge 0}S^j\a_0$ 
be the ring of functions on the corresponding formal group,
and $A$ be a deformation of $A_0$ as a topological Hopf algebra. 
Let $\m$ be the maximal ideal in $A$, i.e. 
the kernel of the projection $A\to k$. 

Let $B$ be a $\d$-copseudotriangular $\O$-structure on $A$.
We define the quasiclassical limit of $B$ 
as follows. 

Let $U_A$ be the h-adic completion
of the direct sum $\oplus_{j\ge 0}h^{-j}\m^j$ (see \cite{Dr1}). 
Then $U_A$ is a quantized 
universal enveloping algebra. 

Let $\a$ be the Lie bialgebra over $k$ which is the quasiclassical limit 
of $U_A$ \cite{Dr1}. Let $\mu,\delta$ be the commutator and the cocommutator
in $\a$, and $\d_0:\a\to\a$ be the quasiclassical limit of $\d$. 
In particular if $\d=D$, then $\d_0=d$
(see \cite{Dr1}, Section 8). 

Consider the pairing $B:\m\o\m\to \O[[h]]$. 
It is clear that $B|_{\m\o\m}=O(h)$, so we define 
$\beta:\m\o\m\to \O$ by $\beta(x,y):=\frac{B(x,y)}{h}\text{ mod } h$. 

Let $\m_0=\m/h\m$. It is clear that $\beta$ descends
to a bilinear form $\beta:\m_0\o\m_0\to \Cal O$, 
which vanishes on $\m_0^2$ on the left 
and on the right. As $\m_0/\m_0^2$ is naturally 
identified with $\a$, we get a bilinear form
$\beta:\a\o\a\to \O$. 

\proclaim{Proposition 1.8} The form $\beta$ satisfies equations (1.13)-(1.15).
\endproclaim

\demo{Proof} Relations (1.13)-(1.15) for $\beta$ 
are easily obtained from relations 
(1.17)-(1.19) for $B$. 
$\square$\enddemo

Thus, Proposition 1.8 states that the quasiclassical limit of $U_A$ is 
endowed with a natural $\d$-copseudotriangular $\O$-structure $\beta$. 
We call $\beta$ the quasiclassical limit of $B$, and $B$ a quantization of 
$\beta$.

Similar definitions and statements apply to weak $\d$-copseudotriangular 
structures.

\subhead 1.7. Factorizations associated to $\d$-copseudotriangular structures
\endsubhead

Consider the field
$F_n=k((u_1))((u_2))...((u_n))$.
For $i<j$, we have a subfield $k((u_i))((u_j))\subset F_n$.
Consider the embedding $k((t))\to k((u))((v))$ 
by the formula
$$
f(t)\to f(u-v):=\sum_{m\ge 0}f^{(m)}(u)(-v)^m/m!.
$$ 
Using this embedding, we define subfields $k((u_i-u_j))\subset k((u_i))((u_j))
\subset F_n$. 

Let $\a$ be a Lie bialgebra 
over $k$ with a weak $\d$-copseudotriangular structure $\beta$. 
Let $\a_{F_n}:=\a\o_kF_n$ be a Lie bialgebra over the field $F_n$ obtained by 
extension of scalars from $\a$. In this section we will define 
a structure of a factored Lie bialgebra on 
the space $\a_\u^n:=\a_{F_n}^{\oplus n}$. 

Let $\a_i=\a_{F_n}$, $i=1,...,n$. 
For any $i<j$, we define a linear map $\a_i\o\a_j\to \a_i\oplus\a_j$ 
by the formula
$$
\mu_{ij}(x,y)=\beta_{13}(\delta(x)\o y)(u_i-u_j)\oplus 
\beta_{12}(x\o\delta(y))(u_i-u_j)\tag 1.23
$$
For $a,b\in\a_\u^n$, such that
$a=\sum_{i=1}^na_i$, $b=\sum_{i=1}^n b_i$,
 $a_i,b_i\in\a_i$, 
set
$$
[a_1+...+a_n,b_1+...+b_n]=\sum_{j=1}^n[a_j,b_j]-\sum_{i<j}
(\mu_{ij}(a_i,b_j)-\mu_{ij}(b_i,a_j)).\tag 1.24
$$

The space $\a_\u^n$ has a natural Lie coalgebra
structure $\delta$, coming form the Lie coalgebra structures 
on $\a_1,...,\a_n$. 

\proclaim{Proposition 1.9} The bracket $[,]$ is a Lie bracket on $\a_\u^n$,
and $\delta$ is a Lie bialgebra structure on $(\a_\u^n,[,])$. 
\endproclaim

\demo{Proof} Skew symmetry of $[,]$ is obvious. The Jacobi identity
and the cocycle condition for $\delta$ is verified by a direct computation.
$\square$\enddemo

Proposition 1.9 shows that any weak $\d$-copseudotriangular structure on $\a$ 
defines a natural structure of a factored Lie bialgebra on $\a_\u^n$. 

Now consider the quantum analogue of this construction.
Let $A$ be a Hopf algebra with a weak copseudoriangular structure $B$. 
Define linear operators
$X_l,X_r,X: A\o A\to A\o A((u))$ by the formula
$$
X_l(a\o b)=B_{13}(\Delta(a)\o\Delta(b)),\
X_r(a\o b)=B_{24}(\Delta(a)\o\Delta(b)), X=X_lX_r^{-1}.\tag 1.25
$$

Because of property (1.20) of $B$, we have $X_l,X_r=1+O(h)$. 

\proclaim{Proposition 1.10} (i) For any $i,j,p,q\in \{1,...,n\}$
$[X_l^{ij}(u),X_r^{pq}(v)]=0$.

Also, $Y=X_r^{-1},X_l$ satisfy the following equations:

(ii) The
quantum Yang-Baxter equation
$$
Y^{12}(u_1-u_2)Y^{13}(u_1-u_3)Y^{23}(u_2-u_3)=
Y^{23}(u_2-u_3)Y^{13}(u_1-u_3)Y^{12}(u_1-u_2).\tag 1.26
$$

(iii) 
$$
Y(a\o bc)=(1\o m)(Y^{12}Y^{13}(a\o b\o c)),\
Y(ab\o c)=(m\o 1)(Y^{23}Y^{13}(a\o b\o c));\tag 1.27
$$

(iv) $Y(1\o a)=1\o a$, $Y(a\o 1)=a\o 1$, $a\in A$.
\endproclaim

\demo{Proof} The identity $[X_l^{ij}(u),X_r^{pq}(v)]=0$
is obvious from the definition. 
The Yang-Baxter equation follows from the properties (1.17)-(1.19)
of $B$.
Identities (iii),(iv) follow directly from properties (1.17),(1.20) of $B$. 
$\square$\enddemo

Let $A_{F_n}$ be the h-adic completion of $A\o_k{F_n}$. $A_i=A_{F_n}$, 
$i=1,...,n$. 
 For $i<j$, let $X_{ij}: A_i\o A_j\to A_i\o A_j$ 
be the operators defined by the formula
$X_{ij}:=X(u_i-u_j)$. 
Proposition 1.10 implies that $X_{ij}$ satisfy properties
(i)-(iii) in Section 1.1, so they define a factored algebra $P_n(A_1,...,A_n,
\{X_{ij}\})$.

\proclaim{Proposition 1.11} The ideal $I(\{X_{ij}\})$ is a Hopf ideal
with respect to the coproduct on $A_1*...*A_n$ induced by the coproduct on 
$A_i$. 
\endproclaim

\demo{Proof} It is convenient to write the
relations of the ideal $I(\{X_{ij}\})$ in the form
$m(X_l(a\o b))=m^{op}(X_r(a\o b))$, $a\in A_i,b\in A_j$, $i<j$. Then it is
easy to verify directly that this relation
is invariant under $\Delta$. 
$\square$\enddemo

Proposition 1.11 shows that 
$P_n(A_1,...,A_n,\{X_{ij}\})$ has a natural Hopf algebra structure. 
We denote this Hopf algebra by $A_\u^n$.

Now let $A$ be an h-formal group, and $U_A$ the corresponding 
QUE algebra. Let $B$ be a $\d$-copseudotriangular structure on $A$. 
Although $B$ does not extend to $U_A$, we have the following proposition.

\proclaim{Proposition 1.12} The operator
$X$ extends to an invertible, h-adically continuous 
operator $X: U_A\o U_A\to U_A\o U_A((u))$.
\endproclaim

\demo{Proof} It is easy to see that for any $x,y\in\m$
$X_l(u)(x\o y)$ and $X_r(u)(x\o y)$ belong to the coset $B(x,y)(u)1\o 1+
\m\o \m$. On the other hand, as $X_l,X_r=1+O(h)$, we have
$X=X_lX_r^{-1}=1+X_l-X_r+O(h^2)$. Therefore, $X(\m\o \m)\subset\m\o\m$. 
Using Proposition 1.10, we see that 
$X(\m^r\o \m^s)\subset \m^r\o \m^s$ for any $r,s\ge 0$. 
Therefore, $X$ extends to $U_A\o U_A$. Similarly, we can extend 
$X^{-1}=X_rX_l^{-1}$ to $U_A\o U_A$. 
$\square$\enddemo

Proposition 1.12 enables us to construct the factored Hopf algebra
$P_n(U_{A_1},...,U_{A_n},\{X_{ij}\})$. 
We denote this Hopf algebra by  
$U_{A_\u^n}$.

{\bf Remark.} It is easy to see that $A_\u^n$ is an h-formal group
over $F_n$, and $U_{A_\u^n}$ is the corresponding QUE algebra, so this notation
is consistent with the previous notation.

Now let us consider copseudotriangular structures and factorizations
which are defined over the ring of functions on some algebraic variety. 

Let $\Sigma$ be a connected
1-dimensional algebraic group over $k$ with a fixed invariant 
differential $du$, 
$\Gamma$ a finite subset of $\Sigma(k)$, 
and $\O$ be the algebra of regular functions on $\Sigma\setminus\Gamma$. 
We can regard $\O$ as a subalgebra in $k((u))$ using the canonical
formal parameter $u$ near the origin (the parameter whose 
differential is $du$). 

Let $\Sigma_n(\Gamma)$ be the variety
of all $\z=(z_1,...,z_n)\in\Sigma^n$ 
such that $z_i-z_j\notin\Gamma$. 
Let $\O_n$ be the ring of regular functions on $\Sigma_n(\Gamma)$.
We have a natural embedding 
$\O_n\to F_n$, which acts by taking the Laurent expansion 
of a function $f\in\O_n$ near the origin, consecutively 
in the variables $z_n,...,z_1$.  

Let $\beta$ be a weak $\d$-copseudotriangular 
$\O$-structure on a Lie bialgebra $\a$. 
Then the Lie bialgebra $\a_\u^n$ over $F_n$ has a
natural $\O_n$-structure. 
Indeed, the $\O_n$-submodule $\a_\u^{\O_n}:=
\oplus_{i=1}^n{\a\o_k\O_n}\subset \a_\u^n$
is a Lie bialgebra over $\O_n$, and $\a_\u^n=\a_\u^{\O_n}\o_{\O_n}F_n$. 
For any $\z\in 
\Sigma_n(\Gamma)(k)$, define 
the Lie bialgebra $\a_\z:=\a_\u^{\O_n}/I(\z)$, where 
$I(\z)\subset \O_n$ is the ideal of functions vanishing at $\z$.
Then $\a_\z$ is a factored Lie bialgebra over $k$, with $n$ factors 
isomorphic to $\a$. 

Similarly, one defines factored Hopf algebras $A_\u^{\O_n}$,
$A_\z$, $U_{A_{\bold u}^{\O_n}}$, $U_{A_\z}$. 

\proclaim{Proposition 1.16}
If $U_A$ is a quantization of a Lie bialgebra
$\a$, and $B$ is a quantization of a weak $\d$-copseudotriangular structure $\beta$
on $\a$, then the QUE algebra $U_{A_\u^n}$ is a quantization
of the Lie bialgebra $\a_\u^n$. 
If in addition $\beta$, $B$ are $\O$-structures, then 
$U_{A_\u^{\O_n}},U_{A_\z}$ are quantizations of $\a_\u^{\O_n},\a_\z$.  
\endproclaim

\demo{Proof} Easy.
$\square$\enddemo

\subhead 
1.8. Quantization of weak $\d$-copseudotriangular structures
\endsubhead

In this section we will show that any
weak $\d$-copseudotriangular 
structure on a Lie bialgebra $\a$ admits a quantization. 

Let $\a$ be a Lie bialgebra over $k$ with a weak $\d$-copseudotriangular structure 
$\beta$. Consider the Lie bialgebra $\a_\u^n$ 
over $F_n$ defined above. 
Define a linear map $\theta_n: \a_\u^n\to \a_{F_n}^*$ by
$$
\theta_n(a_1+...+a_n)(b)=\beta(a_1,b)(u_1)+...+\beta(a_n,b)(u_n).\tag 1.28
$$

\proclaim{Proposition 1.17} 

(i) $\theta_n$ is a homomorphism of Lie bialgebras $\a_\u^n\to
\a_{F_n}^{*op}$
(where ${\frak b}^{op}$ is ${\frak b}$ with opposite cocommutator, 
for any Lie bialgebra ${\frak b}$). 

(ii) $\theta_1(\alpha_v(x))(u)=\alpha_{-v}^*(\theta_1(x)(u))=\theta_1(x)(u+v)$.

(iii) $\theta_n(a_1+...+a_n)(u_1,...,u_n)=\theta_1(a_1)(u_1)+...+
\theta_1(a_n)(u_n)$, $a_i\in \a_i$.

(iv) $\theta_1$ is injective if and only if $\beta$ is left-nondegenerate. 
\endproclaim

\demo{Proof} Properties (i)-(iv) follow from (1.13),(1.15),(1.21). 
We check only the case of property (i),
which is less obvious than others. 
For $n=1$, (i) follows from (1.13). 
Consider the case $n\ge 2$. It is clear that we can assume $n=2$. 
In this case (i) reduces to the identity
$$
\theta_2(\mu(x,y))=[\theta_1(y)(v),\theta_1(x)(u)], 
x,y\in \a.\tag 1.29
$$
Using the definition (1.23) of $\mu$, we rewrite (1.29) in the form
$$
\theta_1(u)(\beta_{13}(\delta(x)\o y)(u-v))+ 
\theta_1(v)(\beta_{12}(x\o\delta(y))(u-v))=[\theta_1(v)(y),\theta_1(u)(x)].
\tag 1.30
$$
Evaluating both sides of (1.30) on an element $z\in\a$, we rewrite
(1.30) in the form
$$
\beta_{13}(u-v)\beta_{24}(u)(\delta(x)\o y\o z)+\beta_{12}(u-v)\beta_{34}
(v)(x\o\delta(y)\o z)+\beta_{13}(u)\beta_{24}(v)(x\o y\o \delta(z))=0.
\tag 1.31
$$
On the other hand,  
and using (1.13), (1.15), we can rewrite (1.21) in the form
$$
\beta_{13}(u)\beta_{24}(u+v)(\delta(x)\o y\o z)+\beta_{12}(u)\beta_{34}
(v)(x\o\delta(y)\o z)+\beta_{13}(u+v)\beta_{24}(v)(x\o y\o \delta(z))=0.
\tag 1.32
$$
It is easy to see that (1.32) is transformed to (1.31) by the change
of variable $u=u'-v',v=v'$. The proposition is proved. 
$\square$\enddemo

Let $\a$ be a Lie bialgebra over $k$, $\d_0$ a derivation of $\a$,
$\beta:\a\o\a\to k((u))$ a weak $\d_0$-copseudotriangular structure on $\a$.
Let $U_h$ be the functor of quantization of Lie bialgebras (see section 1.2).
Let $U_A=U_h(\a)$, $\d=U_h(\d_0)$, and 
$A$ be the h-formal group corresponding to $U_A$.
(\cite{Dr1}, Section 7). 

\proclaim{Theorem 1.18} The Hopf algebra $A$ admits a 
$\d$-copseudotriangular structure $B$, which is a quantization of $\beta$, 
such that $B$ is left-nondegenerate if and only if so is $\beta$.
\endproclaim

\demo{Proof} 
Since $U_h$ is a functor, the Hopf algebra 
$U_h(\a_\u^n)$ admits a factorization $U_{A_1}\o...\o U_{A_n}$, where 
$U_{A_i}$ is the h-adic completion of $U_A\o_kF_n$. 
 
By Proposition 1.17, the linear map  
$\theta_n$ constructed in the previous section,
is a homomorphism of Lie bialgebras. 
Therefore, by functoriality of quantization
 \cite{EK2}, it defines a homomorphism of Hopf algebras 
$\tilde\Theta_n: U_h(\a^n_\u)\to U_h(\a_{F_n}^{*op})$,
by $\tilde\Theta_n=U_h(\theta_n)$.

Therefore, 
Proposition 1.17 implies the following 
properties of $\tilde\Theta_n$. 

(i) $\tilde\Theta_n$ is a homomorphism of Hopf algebras 
$U_h(\a^n_\u)\to U_h(\a_{F_n}^{*op})$. 

(ii) $\tilde\Theta_1(\alpha_v(x))(u)=
\alpha_{-v}^*(\tilde\Theta_1(x)(u))=\Theta_1(x)(u+v)$.

(iii) $\tilde\Theta_n(a_1\o...\o a_n)(u_1,...,u_n)=\tilde\Theta_1(a_1)(u_1)...
\tilde\Theta_1(a_n)(u_n)$, $a_i\in U_{A_i}$.

(iv) $\tilde\Theta_1$ 
is injective if and only if $\beta$ is left-nondegenerate. 

Let $I$ be the kernel of $\tilde\Theta_1$. It is clear that $I$ 
is a Hopf ideal in $U_h(\a)$. 

For any Lie bialgebra $\b$, we have a natural Hopf algebra isomorphism 
$\psi: U_h(\b^{*op})\to U_h(\b)^{*op}$ defined as follows. 
Consider the natural Lie bialgebra maps $\eta:\b\to D(\b)$, 
$\eta_*:\b^{*op}\to D(\b)$, where $D(\b)= \b\oplus \b^{*op}$
is the double of $\b$. 
These maps are given by $\eta(x)=(x,0)$, $\eta_*(y)=
(0,y)$. Let $U_h(\eta): U_h(\b)\to U_h(D(\b))$,
$U_h(\eta_*): U_h(\b^{*op})\to U_h(D(\b))$ be their quantizations. 
Let $\Cal R$ be the universal R-matrix of $U_h(D(\b))$. By \cite{EK1,EK2}, 
we have $\Cal R\in Im U_h(\eta)\o Im U_h(\eta_*)$. Thus,
$\Cal R$ can be regarded as an element of $U_h(\b)\o U_h(\b^{*op})$, 
hence as a linear map $\phi: U_h(\b^{*op})\to U_h(\b)^{*op}$.
It is clear that $\phi$ is an isomorphism of Hopf algebras,    
so we can define $\psi=\phi^{-1}$. 

Now set $\Theta_n=\psi\circ \tilde\Theta_n$. It is clear that 
$\Theta_n$ is a homomorphism of Hopf algebras  
$U_h(\a^n_\u)\to U_h(\a_{F_n})^{*op}$, which satisfies properties 
(ii)-(iv) above.  
 
Define the bilinear form $B:A\o A\to k((u))[[h]]$ by
$B(a,b)=\Theta_1(a)(b)$. 

We claim that $B$ is a weak $\d$-copseudotriangular structure.

Indeed, identity (1.17) follows from the fact that $\Theta_1$ is a 
homomorphism. Property (iv) of $\Theta_n$  
implies (1.19). Property (1.20) is clear.  
It remains to establish property (1.22).

Let us prove (1.22). Let $\a^2_\u$ be as above, and 
$D(\a)=\a\oplus \a^{*op}$ be the double of 
$\a$. Define a linear map 
$\chi: \a^2_\u\to D(\a)\o F_2$ defined by 
\linebreak $(a,b)\to (b,\theta_1(a)(u_1-u_2))$. 

{\bf Lemma 1.19.} $\chi$ is a Lie bialgebra homomorphism. 

{\it Proof.} Straightforward. 

Lemma 1.19 allows us to define a homomorphism 
$\hat\chi=U_h(\chi): U_h(\a^2_\u)\to 
U_h(D(\a))=U_h(\a)\o U_h(\a^{*op})$.  
This homomorphism satisfies the equation
 $\hat\chi(ab)=\tilde\Theta_1(a)b$, where 
$a$ is from $U_h(\a)_1$ and $b$ from $U_h(\a)_2$, and 
$U_h(\a)_{1,2}$ denote the first and second components of 
$U_h(\a)$ in 
$U_h(\a^2_\u)$. This shows that we have:

{\bf Lemma 1.20.} The linear map  
$\xi: U_h(\a^2_\u)\to U_h(\a)\o U_h(\a)^{*op}=D(U_h(\a))$ given by 
$\xi(ab)=\Theta_1(a)b$, where 
$a$ is from $U_h(\a)_1$ and $b$ from $U_h(a)_2$, 
is a Hopf algebra homomorphism. 

Lemma 1.20 allows us to get information about commutation relations 
between the two components of $U_h(\a)$ inside of $U_h(\a^2_\u)$. 

{\bf Lemma 1.21.} Modulo $I\o U_h(\a)$, the multiplication in 
$U_h(\a^2_\u)$ satisfies the equation 
$m^{op}(X_r(u_1-u_2)(a\o b))=m(X_l(u_1-u_2)(a\o b))$, where $X_{l,r}(u):
U_h(\a)^{\o 2}\to U_h(\a)^{\o 2}$ 
are the operators defined by $B$ according to formula (1.25). 

{\bf Remark.} Note that since $I$ is a Hopf ideal, the operators 
$X_{l,r}$ preserve $I\o U_h(\a)$. 

{\it Proof} It is enough to check that the relation of Lemma 1.21 
is satisfied after applying the map $\xi$. This easily follows 
from property (1.17) of $B$ and the 
commutation relation between $U_h(\a)$ and its dual in the double. 

Now we can prove property (1.22) of $B$. For this purpose we will apply  
property (iii) of $\Theta_n$ for $n=2$. Namely, using Lemma 1.21, we see 
that the fact that the map $\Theta_2$ defined by $\Theta_2(ab)(u_1,u_2)=
\Theta_1(a)(u_1)\Theta_1(b)(u_2)$ 
is a homomorphism of Hopf algebras gives exactly (1.22).  

Thus, $B$ is a weak $\d$-copseudotriangular structure on $U_A$. 
It is easy to show that $B$ is a quantization of $\beta$.
It is also clear that $B$ is left-nondegenerate iff $\Theta_1$
is injective,
so the nondegeneracy of $B$ is equivalent to the nondegeneracy of $\beta$. 
The theorem is proved.

{\bf Remark.} Since the Lie algebra $\a$ is allowed to be 
infinite-dimensional (cf. Proposition 1.7), it is necessary 
to clarify what is meant by $U_h(\a^{*op})$ and $U_h(\a)^{*op}$. 
As vector spaces, these algebras are equal to 
$(\oplus_{n\ge 0}(S^n\a)^*)[[h]]$, with the operations defined in the same 
way as in the finite-dimensional case.   
$\square$\enddemo

The form $B$ constructed in the proof of Theorem 1.18 will be
denoted by $U_h(\beta)$.

\proclaim{Proposition 1.22} Suppose $\beta$ is a left-nondegenerate 
$\d_0$-copseudotriangular structure on $\a$. 
Then the factorization of $U_h(\a^n_\u)$ into 
$n$ copies of $U_h(\a)$ is defined by the 
$\d$-copseudotriangular structure $B=U_h(\beta)$. 
\endproclaim

\demo{Proof} This follows from Lemma 1.21 and the fact that $I=0$ in the 
nondegenerate case.  
\enddemo

\proclaim{Proposition 1.23} Any weak copseudotriangular structure
on a Lie bialgebra admits a quantization.
\endproclaim

\demo{Proof} The Proposition follows from Theorem 1.18
and the equality $U_h(d)=D$, which is proved in Appendix A
(Proposition A3). 
$\square$\enddemo

\head 2. Lie bialgebras of functions on a curve with punctures\endhead

In this chapter, we will give some important examples 
of $\d$-copseudotriangular structures, which arise from 
solutions of the classical Yang-Baxter equations with a spectral paramater. 

\subhead 2.1. Classical r-matrices\endsubhead

In this section, we remind the definition
of certain Lie bialgebras associated with
a rational or elliptic curve with punctures, which was introduced by Drinfeld,
\cite{Dr1}. As this material is known, and proofs are easy,
we omit proofs of most statements.

Let $\g$ be a finite-dimensional simple Lie algebra over $k$,
$\<,\>$ be the Killing form on $\g$ divided by $2$.
Let $\Omega\in (S^2\g)^\g$ be the inverse element to 
$\<,\>$.

Let $\bg:=\g((t))$ be the Lie algebra of formal
Laurent series with values in $\g$. The algebra $\bg$ has a natural 
nondegenerate invariant inner product $(,)$ defined by $(a,b)=
Res(\<a(t),b(t)\>)$, where $Res(\sum a_mt^m):=a_{-1}$. 
Denote by $\bg^n$ the direct sum of $n$ copies of $\bg$. 
The inner product on $\bg$ induces one on $\bg^n$. 

Let $\bg_-:=\g[[t]]$, and $\bg_-^n$ be the direct sum of $n$ copies
of $\bg_-$. It is clear that for any $n\ge 1$ $\bg_-^n$ is a 
Lie subalgebra in $\bg^n$, isotropic under the inner product $(,)$.
  
Let $r\in\g\o\g((u))$, $r(u)=\frac{\Omega}{u}+O(1)$. 

\proclaim{Definition} One says that $r(u)$ is
a classical $r$-matrix if it satisfies the classical 
Yang-Baxter equation
$$
[r^{12}(u_1-u_2),r^{13}(u_1-u_3)]+
[r^{12}(u_1-u_2),r^{23}(u_2-u_3)]+
[r^{13}(u_1-u_3),r^{23}(u_2-u_3)]=0.\tag 2.1
$$
\endproclaim

Belavin and Drinfeld \cite{BD} showed that any classical r-matrix
also satisfies the unitarity condition
$$
r(u)=-r^{21}(-u).\tag 2.2
$$

\subhead 2.2. Lie bialgebras associated to classical r-matrices\endsubhead

Given a classical r-matrix $r(u)$ on $\g$, one can construct an infinite 
dimensional Lie algebra $\g(r)$ as follows. We regard $r^{21}(t-u)$ as an 
element of $\g\o\g((t))[[u]]$, expanding $\frac{1}{t-u}$ as
$\sum_{m\ge 0}t^{-m-1}u^m$. Define
$$
\g(r):=\{Res(\<X(u)\o 1,r^{21}(t-u)\>): 
X(u)\in u^{-1}\g[u^{-1}]\}\subset \g((t)).
\tag 2.3
$$
We have the following well known propositions.

\proclaim{Proposition 2.1} $\g(r)$ is a Lie subalgebra of
$\g((t))$ isotropic under $\<,\>$. 
\endproclaim

\demo{Proof} The first statement follows from (2.1), and the second one from 
(2.2).  
$\square$\enddemo

\proclaim{Proposition 2.2} $\bg=\g(r)\oplus \bg_-$, and
$(\bg,\g(r),\bg_-)$ is a Manin triple. 
\endproclaim

\demo{Proof} Clear.
$\square$\enddemo

Thus, $\g(r)$ has a natural Lie bialgebra structure. 

The commutation 
and cocommutation relations of $\g(r)$ have the following convenient explicit
representation. Define the generating function
$\tau(u):=r^{21}(t-u)\in \g\o\g(r)[[u]]$. 

\proclaim{Proposition 2.3} One has
$$
\gather
[\tau^{13}(u),\tau^{23}(v)]=
[r^{12}(u-v),\tau^{13}(u)+\tau^{23}(v)],\\
\delta(\tau(u))=[\tau^{12}(u),\tau^{13}(u)].
\tag 2.4\endgather
$$ 
\endproclaim

{\bf Remark 1.} Here and below we use 
the usual notation  $\tau^{12}:=\tau\o 1$, 
$\tau^{23}:=1\o \tau$, and so on. 

{\bf Remark 2.} It is obvious that 
$[\Omega^{12},\tau^{13}(u)+\tau^{23}(u)]=0$, 
so the right hand side of the first equation in (2.4) 
is regular at $u,v=0$, and thus (2.4) really defines a bracket on
$\g(r)$.

\demo{Proof}
Relations (2.4) follow from the classical Yang-Baxter equation for $r(u)$. 
$\square$\enddemo

\subhead 2.3. The canonical derivation\endsubhead

Let $f\in\g\o\g$ be the free term in the Laurent expansion 
of $r$, i.e. $r(u)=\frac{\Omega}{u}+f+O(u)$.
Let $\mu$ be the commutator in $\g$, and 
$\rho_r=\frac{1}{2}\mu(f)$.

\proclaim{Proposition 2.4} The canonical derivation $d$ on 
$\g(r)$ (see Section 1.3) equals $-\frac{d}{dt}+\ad(\rho_r)$. 
\endproclaim

\demo{Proof} Since $\g(r)$ is a subbialgebra of $\g((t))$, it is
enough to show that $d=-\frac{d}{dt}+\ad(\rho_r)$ in $\g((t))$. 

The cobracket in
$\g((t))$ is given by 
$$
\gather
\delta(xu^n)=[xu^n\o 1+1\o xv^n,r(u-v)]=\delta_0(xu^n)+\phi(xu^n),
\\ \delta_0(xu^n):=[xu^n\o 1+1\o xv^n,\frac{\Omega}{u-v}], \ \phi(xu^n):=
[xu^n\o 1+1\o xv^n,f+O(u-v)].\tag 2.5\endgather
$$
Let $d_0=
-\frac{1}{2}\mu\circ \delta_0$. Then from (2.5) we obtain
$$
d=d_0+\ad(\rho_r).\tag 2.6
$$

It remains to show that $d_0=-\frac{d}{dt}$.
We have 
$$
d_0(xu^n)=-\frac{1}{2}\mu([x\o 1,\Omega]\frac{u^n-v^n}{u-v}).\tag 2.7
$$
The normalization of $\Omega$ is such that
$\mu([x\o 1,\Omega])=2x$. Therefore, the right hand side
of (2.7) equals $-nxu^{n-1}$, as desired.
$\square$\enddemo

\proclaim{Proposition 2.5} Any classical r-matrix $r(u)$ is invariant
under the adjoint action of $\rho_r$, and the function
$\tilde r(u):=(e^{-u\rho_r}\o 1)r(u)(e^{u\rho_r}\o 1)$ is also a
classical $r$-matrix.
\endproclaim

\demo{Proof} It is clear that $(d\o 1+1\o d)(r(u-v))=0$, as $d$ preserves
the structure of a Manin triple in $\g((t))$. By Proposition 2.4, this imples 
that $[\rho_r\o 1+1\o\rho_r,r(u)]=0$, as desired. 

Because of this, we can write
$$
\tilde r(u-v)=(\Ad e^{-u\rho_r}\o \Ad e^{-v\rho_r})(r(u-v)).\tag 2.8
$$
It is clear from this formula that $\tilde r$ is a classical r-matrix. 
$\square$\enddemo

Proposition 2.5 shows that the operator
$\ad(\rho_r)$ is a derivation of the Manin triple $(\bg,\g(r),\bg_-)$. 
Define $\d=d-\ad(\rho_r)$. This is also a derivation of this Manin triple. 
We have $\d=-\frac{d}{dt}$, i.e. $\d(\tau(u))=\frac{d\tau(u)}{du}$.

\subhead 2.4. $\d$-copseudotriangular structure on $\g(r)$\endsubhead
 
It turns out that the Lie bialgebra $\g(r)$ has a natural 
$\d$-copseudotriangular structure. Namely, set
$$
\beta_{3}(\tau^{13}(v),\tau^{23}(w))(u)=-r^{12}(v-w+u).\tag 2.9
$$

\proclaim{Proposition 2.6} $\beta$ is a $\d$-copseudotriangular
structure on $\g(r)$.
\endproclaim

\demo{Proof} 
Formulas (1.13) follow from (2.4) and the Yang-Baxter equation for $r$.
Property (1.15) follows from Proposition 2.4. Finally, formula
(1.14) follows from the identity
$\delta(x(u))=[x(u)\o 1+1\o x(v),r(u-v)]$. 
$\square$\enddemo

\subhead 2.5.  Factored Lie bialgebras associated with $\g(r)$\endsubhead

Let $\Sigma$ be a 1-dimensional connected algebraic group over $k$
(i.e. $\Bbb G_a$, $\Bbb G_m$, or an elliptic curve), 
and $du$ be an invariant differential on $\Sigma$.
The differential $du$ defines a canonical formal parameter near any point
of $\Sigma$. Let $u$ be the corresponding 
formal parameter at the origin.

The following deep theorem is due to Belavin and Drinfeld \cite{BD}. 

\proclaim{Theorem 2.7} Let $r$ be a classical $r$-matrix.
Then there exists a unique 1-dimensional algebraic group $\Sigma$
such that the formal series $r(u)$ is the Laurent expansion with respect to 
$u$ of a rational function on $\Sigma$ with values in
$\g\o\g$, and the stabilizer of $r(u)$ in $\Sigma$ is trivial. 
\endproclaim

Let $\Gamma$ be the set of poles of $r(u)$ in $\Sigma$.
Theorem 2.7 shows that the $\d$-copseudotriangular
structure $\beta$ on $\g(r)$ is in fact an $\O$-structure, where
$\O=\O(\Sigma\setminus\Gamma)\subset k((u))$. 

Let
$\z=(z_1,...,z_n)\in\Sigma(k)^n$ be such that $z_i-z_j\notin\Gamma$ for
$i\ne j$. 
The form $\beta$ allows us to define the factored Lie bialgebras 
$\g(r)_\u^{\O_n}$ and $\g(r)_\z$, as in Section 1.7. 

It is easy to write down explicitly the commutator in $\g(r)_\z$.
Recall that $\g(r)_\z$ is equal to $\g(r)^{\oplus n}$ 
as a vector space. 
Let $\tau_i(u)$ be the series $\tau(u)$ for the $i$-th summand $\g(r)$
in this direct sum. 

\proclaim{Proposition 2.8} The bracket in $\g(r)_\z$ has the form
$$
[\tau^{13}_i(u),\tau^{23}_j(v)]=
[r^{12}(u-v+z_i-z_j),\tau^{13}_i(u)+\tau^{23}_j(v)].
\tag 2.10
$$ 
\endproclaim

\demo{Proof} Formula (2.10) is obtained by substitution of (2.4) in (1.23).
$\square$\enddemo

The formula for the cobracket in $\g(r)_\z$ 
follows directly from the definition:
$$
\delta(\tau_i(u))=[\tau_i^{12}(u),\tau_i^{13}(u)].
\tag 2.11
$$ 

\subhead 2.6. The Manin triple associated to $\g(r)_\z$\endsubhead

Consider the map $f_\z:\g(r)_\z\to \bg_n$ defined by
$$
f_\z(\tau_i(u))=(r(t-u+z_1-z_i),...,r(t-u+z_n-z_i)).\tag 2.12
$$
It is easy to check that
 $f_\z$ is a Lie algebra embedding.
Thus we can regard $\g(r)_\z$ as a Lie subalgebra in $\bg^n$. 

\proclaim{Proposition 2.9} (i) The algebra $\g(r)_\z$ is isotropic in $\bg^n$
with respect to the form $(,)$. 

(ii) $\bg^n=\g(r)_\z\oplus \bg_-^n$.

(iii) The form $(,)$ defines a linear isomorphism 
$\bg_-^n\to \g(r)_\z^*$.
\endproclaim

\demo{Proof} (i) follows from the unitarity property of $r(u)$.
(ii) is clear. (iii) follows from (ii). 
$\square$\enddemo

\proclaim{Corollary 2.10} The triple of Lie algebras
$(\bg^n,\g_\z,\bg_-^n)$
is a Manin triple.
\endproclaim

\proclaim{Proposition 2.11} The Lie bialgebra structure on $\g(r)_\z$
coming from the Manin triple of Corollary 2.10 coincides with the one
defined by formula (2.11).
\endproclaim

\subhead 2.7. Examples\endsubhead

{\it Example 1: Yang's rational r-matrix.} 

Let $r(u)$ be the Yang's r-matrix $r(u)=\frac{\Omega}{u}$.
In this case, $\Sigma=\Bbb G_a$, and $u$ is an affine coordinate.
The Lie algebra $\g(r)$ coincides with
the Lie algebra $\a=t^{-1}\g[t^{-1}]$. As a Lie bialgebra, $\g(r)$ is 
the graded dual to the Yangian Lie bialgebra, defined by Drinfeld
in \cite{Dr1}, with opposite commutator. 

Let $\z=(z_1,\dots,z_n)\in k^n$, 
$z_i\ne z_j$ for $i\ne j$. 
Let $R_\z$ be the algebra of all rational functions in one variable 
over $k$
with no poles outside of the points $z_1,\dots,z_n$ which vanish at 
infinity. Consider the Lie algebra $\a_\z=\g\o R_\z$, and the Lie algebra 
embedding $\phi_\z: \a_\z\to \bg^n$, defined by the formula
$\phi_\z(f)=(\phi_{z_1}(f),...,\phi_{z_n}(f))$, where 
$\phi_{z_i}(f)\in\bg$ is the Laurent expansion of $f$ at $z_i$.   
Using this embedding, we can regard $\a_\z$ as a Lie subalgebra
in $\bg^n$.

\proclaim{Proposition 2.12} $\g(r)_\z=\a_\z$. 
\endproclaim

Let $\tilde r(u)$ be any classical r-matrix. Then
$\lim_{\e\to 0}\e\tilde r(\e u)$ equals the Yang's r-matrix $r(u)$.
Therefore, the bialgebras $\g(r)$, $\g(r)_\z$ can be obtained as
a degeneration of $\g(\tilde r)$, $\g(\tilde r)_\z$. 

{\it Example 2: Trigonometric r-matrices.}

Suppose that the Lie algebra
$\g$ is endowed with a standard
polarization $\g=\n^-\oplus\hh\oplus\n^+$, where
$\hh$ is a Cartan subalgebra, $\n^\pm$ the nilpotent subalgebras.

Let $L=L_++\frac{1}{2}L_0\in\g\o\g$,
where $L_+$ is the canonical element in $\n^+\o \n^-$, and
$L_0\in\hh\o\hh$ is the element dual to $\<,\>|_{\hh}$. 
It is easy to check that the element $L$ is a constant solution to the
classical Yang-Baxter equation. Therefore, the element
$$
r(u)=\frac{L^{21}e^u+L}{e^u-1},\tag 2.13
$$ 
is a classical r-matrix. 
In this case,
$\Sigma=\Bbb G_m=\Bbb A^1\setminus 0$, $u=\ln x$,
where $x$ is the affine coordinate
on $\Bbb A^1$. 

Let $\a$ be the Lie algebra of $\g$-valued rational functions $X(x)$ 
on $\Sigma$,
such that 

(i) $X$ has no poles on $\Sigma$ outside of the identity;

(ii) $X$ is regular at the infinite points $0,\infty$ of $\Sigma$;

(iii) $X(0)\in\hh\oplus \n_-$, $X(\infty)\in\hh\oplus\n_+$, and 
$X(0)+X(\infty)\in\n_+\oplus\n_-$. 

Consider the embedding $\phi: \a\to \bg$, which assigns to
$X\in\a$ its Laurent expansion at the identity with respect to $u$.
This embedding allows us to regard $\a$ as
a subalgebra in $\bg$. 

\proclaim{Proposition 2.13} $\g(r)=\a$. 
\endproclaim

Analogously to the construction of $\g(r)$, 
for any $\z\in\Sigma^n$, $z_i-z_j\ne 0$, $i\ne j$, define 
$\a_\z$ to be the Lie algebra of $\g$-valued rational functions $X(x)$ 
on $\Sigma$,
such that 

(i) $X$ has no poles on $\Sigma$ outside of $z_1,...,z_n$;

(ii) $X$ is regular at the infinite points $0,\infty$ of $\Sigma$;

(iii) $X(0)\in\hh\oplus \n_-$, $X(\infty)\in\hh\oplus\n_+$, and 
$X(0)+X(\infty)\in\n_+\oplus\n_-$. 

Consider the embedding $\phi_\z: \a_\z\to \bg^n$, which assigns to
$X\in\a_\z$ the collection of
its Laurent expansions at $z_i$ with respect to the paramaters
 $u_i=\ln (x/x(z_i))$. This embedding allows us to regard $\a_\z$ as
a subalgebra in $\bg^n$. 

\proclaim{Proposition 2.14} $\g(r)_\z=\a_\z$. 
\endproclaim

{\it Example 3: Elliptic r-matrices.} 

Let $\g=\frak{sl}_N$, $\<X,Y\>=N\text{tr}(XY)$.
Let $A:\g\to\g$ be the conjugation by the cyclic permutation
$(12...N)$, and $B:\g\to\g$ be the conjugation by the matrix
$\text{diag}(1,\e,...,\e^{N-1})$, where $\e$ is a primitive 
$N$-th root of $1$.   
Let $H$ be the subgroup in $Aut(\g)$ generated by $A,B$.
This group is isomorphic to $\Z/N\Z\times \Z/N\Z$.

Let $\Sigma$ be an elliptic curve over $k$. 
Let $\Gamma$ be the group 
of points of order $N$ in $\Sigma$. Fix an isomorphism
$\l: \Gamma\to H$. Let $r$ be a $\g$-valued rational function on $\Sigma$
of the form $r(u)=\frac{\Omega}{u}+O(1)$, such that
$r(z+\gamma)=(1\o\l(\gamma))(r(z))$, $\gamma\in\Gamma$, and $r$ is regular 
outside of $\Gamma$.
It is easy to see that such a function is unique. 
It is known \cite{BD} that it is a classical $r$-matrix.

Let $\a$ be the Lie algebra of all $\g$-valued rational functions $X$
on $\Sigma$ such that $X$ has no poles outside of $\Gamma$, and 
$X(z+\gamma)=\l(\gamma)(X(z))$, $\gamma\in\Gamma$. 
Define a Lie algebra embedding
$\phi:\a\to\bg$ by assigning to a function $X$ its Laurent expansion at the 
origin. We identify $\a$ with its image under this embedding. Then we have

\proclaim{Proposition 2.15} $\g(r)=\a$. 
\endproclaim

Analogously, for any $\z\in\Sigma^n$, $z_i-z_j\notin\Gamma$, $i\ne j$, define 
$\a_\z$ to be the Lie algebra of $\g$-valued rational functions $X$ 
on $\Sigma$,
such that $X$ has no poles outside of $\cup_i(z_i+\Gamma)$, and 
$X(z+\gamma)=\l(\gamma)(X(z))$, $\gamma\in\Gamma$.

Consider the embedding $\phi_\z: \a_\z\to \bg^n$, which assigns to
$X\in\a_\z$ the collection of
its Laurent expansions at $z_i$. 
This embedding allows us to regard $\a_\z$ as
a subalgebra in $\bg^n$. 

\proclaim{Proposition 2.16} $\g(r)_\z=\a_\z$. 
\endproclaim

\subhead 2.8. Poisson groups\endsubhead

In this section we consider Poisson groups corresponding to the constructed 
Lie bialgebras. They will be useful in the next Chapter. 

Assume that $\g$ is the Lie algebra of an affine algebraic group $G$
defined over $k$. 
In this case, define the proalgebraic group $\bG_-:=G[[t]]$.
Any classical $r$-matrix $r(u)$ on $\Sigma$
defines a natural
Poisson-Lie structure on this group, coming from the Manin triple 
of Proposition 2.2. We denote the obtained Poisson group by $G(r)$. 
 
Let $\bG_-^n:=\prod_{i=1}^n\bG_-$. The Lie algebra
of this group is $\bg_-^n$.
Therefore, for any 
element $\z=(z_1,\dots,z_n)\in\Sigma(k)^n$, $z_i-z_j\notin \Gamma$,
we can define a natural 
structure on this group, coming from the Manin triple 
of Corollary 2.10. We denote the obtained Poisson group by $G(r)_\z$. 
For example, if $n=1$ then $G(r)_\z$ coincides with $G(r)$. 
 
At the level of formal groups, it is easy to describe the Poisson
bracket on $G(r)_\z$ explicitly.
Let $F_{\z,r}(G)$ 
be the algebra of regular functions on the Poisson group 
$G(r)_\z$.
Let $\pF_{\z,r}(G)$ be the completion of the algebra
$F_{\z,r}(G)$ at the identity. 
Using the exponential map, we can identify 
$\pF_{\z,r}(G)$ with the algebra topologically
generated by linear functions on the Lie algebra 
$\g_-^n$. So, it is enough to write down
the Poisson bracket of such linear functions. 

Using the dualities of Proposition 2.9 (iii),
we can regard the expression $\<X\o 1,\tau_i(u)\>$, $X\in\g$,
 as a formal series 
whose coefficients are linear functions on $\g_-^n$:
$$
\tau_i(u)(a_1,\dots,a_n)=\<r^{12}(t-u),a_i^2(t)\>,\tag 2.14
$$
Denote by $t_i(u)$ the formal series of linear functions
corresponding to $\tau_i(u)$. 
Then the commutation relations
for $\tau_i(u)$ can be rewritten as Poisson commutation relations
for $t_i(u)$:
$$
\{t^{13}_i(u),t^{23}_j(v)\}=
[r^{12}(u-v+z_i-z_j),t^{13}_i(u)+t^{23}_j(v)].\tag 2.15
$$ 
These relations describe the Poisson bracket on the formal version of the 
group $G(r)_\z$.

\head 3. Quantum groups associated to curves with punctures\endhead

\subhead 3.1.  Copseudotriangular structure of dual Yangians and 
associated factored Hopf algebras\endsubhead

In this section we will explicitly 
quantize the Lie bialgebra $\g(r)_\z$, where 
$r$ is the Yang's r-matrix. 

Let $\g$ be a simple finite-dimensional Lie algebra, $\<,\>$ 
be the Killing form on $\g$ divided by $2$. 
Let $Y(\g)$ be the Yangian of $\g$, defined
by Drinfeld, \cite{Dr2}. It is known \cite{Dr2} that $Y(\g)$ is a quantization
of the Lie algebra $\g[t]$, with the cobracket coming from the Manin triple
$(\bg,\g(r),\bg_-)$, via the natural embedding $\g[t]\to \bg$,
where $r(u)$ is the Yang's r-matrix
$\Omega/u$ (See Chapter 2). Thus, as a $k[[h]]$-module, $Y(\g)$ 
is identified with $U(\g[t])[[h]]$.  

Recall
that $Y(\g)$ has a natural $\Z_+$-grading, with $\deg(Xz^j)=j$, $\deg(h)=1$.  
Let $f: Y(\g)\to k[[h]]$ be a continuous
$k[[h]]$-linear functional. We say that $f$ 
is tempered if it can be represented in the form $\sum_{m\ge 0}f_m$, 
where $f_m$ is homogeneous and 
$\lim_{m\to\infty}\deg f_m=+\infty$. Let $Y(\g)'$ denote
the space of all tempered linear functionals on $Y(\g)$.
The $k[[h]]$-module $Y(\g)'$ has a structure of an associative algebra
dual to the coalgebra structure in $Y(\g)$. Let $\frak m$ be
the ideal in $Y(\g)'$ consisting of all functionals $f$ such that 
$f(1)=O(h)$. Denote by $Y(\g)^*$ the $h$-adic completion of
the direct sum $\oplus_{s\ge 0} h^{-s} {\frak m}^s$. 
It is easy to check that $Y(\g)^*$ inherits
the associative algebra structure from $Y(\g)'$. 
Besides, we have

\proclaim{Lemma 3.1} 

(i) the map $m^*$ dual to the product $m :Y(\g)\o Y(\g)\to Y(\g)$ 
takes $Y(\g)^*$ into $Y(\g)^*\o Y(\g)^*$, and defines a coassociative  
coproduct on the topological algebra $Y(\g)^*$, which is cocommutative mod h;

(ii) the algebra $Y(\g)^*$, equipped with $m^*$, is a quantized
universal enveloping algebra, whose quasiclassical limit is $\g(r)$,
with opposite cocommutator. 
\endproclaim

\demo{Proof} Part (i) is easy.
Part (ii) follows from the fact that the quantized
universal enveloping algebra $Y(\g)$ is a quantization of the Lie
bialgebra $\g[t]$, which is the graded dual to $\g(r)$, 
with opposite coproduct.  
$\square$\enddemo

Denote by $Y^*(\g)$ the Hopf algebra $Y(\g)^*$, with 
opposite product. Then $Y^*(\g)$ is a quantization of $\g(r)$, and
the subalgebra $Y(\g)_{op}^{\prime}\subset Y^*(\g)$ 
(i.e. $Y(\g)'$ with opposite product) is the corresponding h-formal group. 

Let $\Sigma=\Bbb G_a$, $\Gamma=\{0\}$. Then $\O=k[u,u^{-1}]$. Now,
using the pseudotriangular structure on $Y(\g)$ defined by Drinfeld
\cite{Dr2}, we will 
define a copseudotriangular $\O$-structure on $Y(\g)_{op}^{\prime}$. 

\proclaim{Proposition 3.2} \cite{Dr2} There exists a unique series
$$
\RR(u)=1+\sum_{m\ge 1}\RR_mu^{-m}, \RR_m\in Y(\g)\o Y(\g)\tag 3.1
$$
which satisfies the hexagon relations 
$$
(\Delta\o 1)(\RR)=\RR^{13}\RR^{23},
(1\o\Delta)(\RR)=\RR^{13}\RR^{12},\tag 3.2
$$
and the property 
$$
\RR(u)(\alpha_u\o 1)(\Delta(x))=(\alpha_u\o 1)(\Delta^{op}(x))\RR(u),\tag 3.3
$$
where $\alpha_u=e^{uD}$.

This series also satisfies the equations $(\e\o 1)(\RR)=(1\o\e)(\RR)=1$,
$(\alpha_v\o 1)(\RR(u))=(1\o\alpha_{-v})(\RR(u))=\RR(u+v):=
\sum \RR^{(m)}(u)v^m/m!$. 
\endproclaim

The element $\RR$ is called 
the pseudotriangular structure on the Yangian $Y(\g)$.

Define the form $B: Y(\g)_{op}^{\prime}
\o Y(\g)_{op}^{\prime}\to \O[[h]]$ by the formula
$$
B(x,y)(u)=(x\o y)(\RR^{21}(-u))=(x\o y)(\RR(u)^{-1}).\tag 3.4
$$

\proclaim{Proposition 3.3}
The form $B$ is a copseudotriangular $\O$-structure on $Y(\g)_{op}^{\prime}$. 
\endproclaim

\demo{Proof} Properties (1.17)-(1.20) of $B$ are dual to the identities 
of Proposition 3.2. 
$\square$\enddemo

{\bf Remark.} This construction explains the terminology 
``copseudotriangular structure''.

Let $\g[t]$ be the Yangian Lie bialgebra, 
considered by Drinfeld in \cite{Dr1}, Example 3.3. The coproduct 
in this Lie bialgebra is defined by
 $$\delta(x(u))=
[x(u)\o 1+1\o x(v),\frac{\Omega}{u-v}].$$

\proclaim{Proposition 3.4} $Y(\g)=U_h(\g[t])$. 
\endproclaim

\demo{Proof} Both $Y(\g)$ and $U_h(\g[t])$ are graded quantizations
of $\g[t]$. Drinfeld showed that a graded quantization of 
$\g[t]$ is unique up to an isomorphism (see \cite{Dr2}).
$\square$\enddemo

Let $r$ be the Yang's r-matrix $\Omega/u$,
$\g(r)$ be the Lie bialgebra from Example 1 of Section 2.7, and
$\beta$ be the copseudotriangular structure on $\g(r)$ defined in
Section 2.4. 

\proclaim{Proposition 3.5} $Y^*(\g)=U_h(\g(r))$, and $B=U_h(\beta)$.
\endproclaim

\demo{Proof} 
By the definition, $Y^*(\g)$ is the graded dual algebra to 
$Y(\g)$, with opposite product, and $\g(r)$ is the graded 
dual bialgebra to $\g[t]$, with opposite commutator.
Therefore, Proposition 3.4 implies 
that $Y^*(\g)=U_h(\g(r))$. 

Let $B'=U_h(\beta)$. Then $B'$ defines an element $\RR'\in Y(\g)\o Y(\g)((u))$,
such that $B'(x,y)=(x\o y,(\RR^{'21})^{-1})$. We know that $B'$ is a weak 
copseudotriangular structure. Also, since $\beta$ is left-nondegenerate,
so is $B'$. Therefore, $B'$ is a 
copseudotriangular structure. Thus, the element $\RR'$  
satisfies properties (3.2),(3.3). Therefore, by 
Proposition 3.2, $\RR'=\RR$, and thus $B'=B$.
$\square$\enddemo

Let $(z_1,...,z_n)\in k^n$ be a collection of distinct points,
and $Y(\g)_{op}^{\prime}(\z), Y^*(\g)(\z)$ 
be the factored Hopf algebras $A_\z$, $U_{A_\z}$
obtained 
from $A:=Y(\g)_{op}^{\prime},U_A:=Y^*(\g)$
 with the help of the copseudotriangular structure $B$,
as described in Section 1.7. Let $\g(r)_\z$ be the 
factored Lie bialgebra
from Example 1 of Section 2.7.

\proclaim{Proposition 3.6} $Y^*(\g)(\z)=U_h(\g(r)_\z)$. 
\endproclaim

\demo{Proof} Since $\g(r)_\z$ is obtained from $\g(r)$ with the help
of the form $\beta$, by Proposition 1.22 $U_h(\g(r)_\z)$ is obtained from
$U_h(\g(r))$ with the help of the form $U_h(\beta)$. 
Therefore, the statement follows from Proposition 3.5.
$\square$\enddemo

{\bf Remark.} Proposition 3.6 gives an ``explicit'' description
of the quantization $U_h(\g(r)_\z)$, in the sense that it does
not use a Drinfeld associator, which is used in the construction
of $U_h(\a)$ for a general Lie bialgebra $\a$, as in \cite{EK1}.

\subhead 3.2.  R-matrix Hopf algebras\endsubhead

It this section we will consider quantum groups defined 
by FRT-type relations \cite{FRT,RS}.

Let $R(u)\in\End(k^N\o k^N)((u))[[h]]$ be an element
such that $R=1+O(h)$. Following
\cite{FRT,RS}, one associates to $R$
the h-adic completion $F(R)$ of the 
algebra over $k[[h]]$ whose generators are the entries
of formal series $T(u)^{\pm 1}\in \End(k^N)\o F(R)[[u]]$,
$T(u)=T_0+T_1u+...$,
and the defining relations are
$$
\gather
T(u)T(u)^{-1}=T(u)^{-1}T(u)=1,\\
R^{12}(u-v)T^{13}(u)T^{23}(v)=T^{23}(v)T^{13}(u)R^{12}(u-v)
\tag 3.5\endgather
$$

\proclaim{Proposition 3.7} There exists a unique Hopf algebra structure
on $F(R)$ such that 
$$
\Delta(T(u))=T^{12}(u)T^{13}(u).\tag 3.6
$$
\endproclaim

\demo{Proof} Introduce the coproduct on the free algebra generated 
by $T(u)^{\pm 1}$ modulo the first relation of (3.5). 
It is easy to check that the ideal in this algebra
generated by the second relation in (3.5) is a Hopf ideal. 
The proposition is proved. 
$\square$\enddemo

\proclaim{Proposition 3.8} The map $\phi_0: F(R)/hF(R)\to
k[GL_N[[t]]]$ defined by $\phi_0(T(u))(g(t))=g(u)$,
$g(*)\in GL_N(k)[[t]]$,
is a Hopf algebra isomorphism.
\endproclaim
 
\demo{Proof} Clear.
$\square$\enddemo

It is interesting to determine when $F(R)$ is a quantization
of the proalgebraic group $GL_N[[t]]$, i.e. when $F(R)$ is
a flat deformation of the function algebra $k[GL_N[[t]]]$.
Flatness is equivalent to the property that the operator of multiplication by $h$ on 
$F(R)$ is injective. 

In general, the algebra $F(R)$ is not a flat deformation of
the function algebra $k[GL_N[[t]]]$. The following proposition 
gives necessary and sufficient conditions of flatness for a large class of 
$R$-matrices.  

\proclaim{Proposition 3.9} 
(a) If $F(R)$ is a flat deformation of $k[GL_N[[t]]]$
then 
$R(u)$ is a solution of the quantum
Yang-Baxter equation
$$
R^{12}(u_1-u_2)R^{13}(u_1-u_3)R^{23}(u_2-u_3)=
R^{23}(u_2-u_3)R^{13}(u_1-u_3)R^{12}(u_1-u_2),\tag 3.7
$$
and $R^{21}(-u)R(u)$ is a scalar operator of the form $f(u)\Id$, where
$f(u)\in 1+hk((u))[[h]]$, 
such that 
$$
f(u)=f(-u).\tag 3.8
$$ 

(b) Let $R(u,h)\in k[[h,u,h/u]]$, and suppose that 
$\lim_{s\to 0}R(su,sh)=R_{rat}(u)$, where 
$R_{rat}(u)=1-\frac{h(\sigma-1/N)}{N(u-h/N^2)}$
is the rational $R$-matrix. 
Then, if $R$ satisfies (3.7), 
$F(R)$ is a flat deformation of $k[GL_N[[t]]]$.
\endproclaim

\demo{Proof}
The proof is given in Section 3.4.
$\square$\enddemo

{\bf Remark.} Note that the converse to (a) does not hold. 
For example, if $R=1-h\Omega/u$, where $\Omega$ is the Casimir
element of $O(N)\subset GL(N)$, then the conclusion of (a) holds but 
$F(R)$ is not flat. 
This is easily seen from considering the quasiclassical limit 
(one does not get a Poisson-Lie structure on 
$GL_N[[t]]$).   

It is clear that rescaling of $R$, i.e. 
multiplication of $R$ by a scalar in $k((u))[[h]]$,
does not change
the algebra $F(R)$. Therefore, Proposition 3.9 implies that
if $F(R)$ is a flat deformation, 
one can arrange that $R$ satisfies the unitarity condition
$$
R(u)R^{21}(-u)=1.\tag 3.9
$$

We know that if $R(u)$ satisfies 
(3.7) and (3.9) then $R(u)=1-hr(u)+O(h^2)$, where 
$r(u)\in \frak{gl}_N\o\frak{gl}_N((u))$ is a unitary solution of the classical 
Yang-Baxter equation. In addition, if 
the assumption of Proposition 3.9(b) holds, 
we have $r(u)=(\sigma-1/N)/u+O(1)$,
where $\sigma$ is the permutation.
It is explained in Section 2.8 that $r$ defines a 
Poisson group $G(r)$ ($G=GL_N$). Proposition 3.9(b) implies that 
the algebra $F(R)$ is always a quantization of the Poisson group $G(r)$.

Assume that $F(R)$ is a flat deformation. 

Let $A(R)$ be the completion of $F(R)$ with respect to the 
kernel of the counit. Then $A(R)$ is an h-formal group, which
is a flat deformation of the formal group corresponding 
to $GL_N[[t]]$. 
Denote the corresponding QUE algebra by $U(R)$. 

Let $R(u)=1-hr(u)+O(h^2)$, where $r\in \End(k^N\o k^N)((u))$. 
According to Proposition 3.9, $r(u)$ satisfies the classical 
Yang-Baxter equation (2.1), and $r(u)+r^{21}(-u)$ is a scalar. 
Define a Lie bialgebra 
$\frak{gl}_N(r)$ by relations (2.4). 

It is easy to describe the algebra $U(R)$ by generators and relations.
 
Namely, let $t(u)\in \End(k^N)\o U(R)[[u]]$ be defined by the formula
$T(u)=1+ht(u)$. Let $r_*(u)=h^{-1}(1-R(u))$.
It is easy to see that relations 
(3.5),(3.6) are equivalent to the following commutation 
and cocommutation relations:
$$
\gather
[t^{13}(u),t^{23}(v)]=
[r_*^{12}(u-v),t^{13}(u)+t^{23}(v)]+\\
h(r_*^{12}(u-v)t^{13}(u)t^{23}(v)-
t^{23}(v)t^{13}(u)r_*^{12}(u-v)),\\
\Delta(t(u))=t^{12}(u)+t^{13}(u)+
ht^{12}(u)t^{13}(u).\tag 3.10\endgather
$$ 
It is easy to see that 
the QUE algebra $U(R)$ is a quantization of $\frak{gl}_N(r)$,
because relations (3.10) coincide with (2.4) modulo $h$. 

{\bf Remark.} R-matrix algebras were first considered
by Cherednik \cite{Ch1}.

\subhead 3.3. R-matrix realization of dual Yangians\endsubhead

In this section we will describe dual Yangians 
for $\g=\frak{gl}_N$ and $\frak{sl}_N$ in the R-matrix language.

Let $R(u)$ be the Yang's quantum R-matrix $R_Y(u):=
1-\frac{h\sigma}{Nu}$.

It is well known (e.g. see \cite{Ch1}) that $F(R)$ is a flat 
deformation of $k[GL_N[[t]]]$. 

\proclaim{Proposition 3.11} The algebra $U(R_Y)$ is a quantization 
of the Lie bialgebra $\frak{gl}_N(r)$, where $r=\sigma/Nu$ is the Yang's
r-matrix for $\frak{gl}_N$. 
\endproclaim

\demo{Proof} Clear.
$\square$\enddemo

Now we will formulate the analogs of these results for the Lie
algebra $\frak{sl}_N$.

Let $A$ be an algebra over $k$. 
For any matrix $X(u)\in Mat_N(k)\o A[[u^{-1}]]$, define
its quantum determinant by 
$$
\qdet(X)=\sum_{\sigma}(-1)^\sigma 
X(u-\frac{h(N-1)}{2N})_{1\sigma(1)}
X(u-\frac{h(N-3)}{2N})_{2\sigma(2)}\dots X(u+\frac{h(N-1)}{2N})_{N\sigma(N)}.
\tag 3.11
$$

Let $F_0(R_Y)$ denote the quotient of $F(R_Y)$ by the relation
$$
\qdet(T(u))=1,\tag 3.12
$$
It is easy to see that $F_0(R_Y)/hF_0(R_Y)=k[SL_N[[t]]]$.

From the flatness of $F(R_Y)$ it can be deduced that 
$F_0(R_Y)$ is a flat deformation of the function algebra 
$k[SL_N[[t]]]$. 

Let $A_0(R_Y)$ be the h-formal group
obtained by completion of $F_0(R_Y)$ with respect to the kernel
of the counit, and $U_0(R_Y)$ be the QUE algebra corresponding to
$A_0(R_Y)$. Analogously to Proposition 3.11 we have

\proclaim{Proposition 3.12} The algebra $U_0(R_Y)$ is a quantization 
of the Lie bialgebra $\g(r)$, where $\g=\frak{sl}_N$,
and $r=\Omega/u$ is the Yang's
r-matrix for $\g$. 
\endproclaim

\proclaim{Proposition 3.13} $U_0(R_Y)$ is isomorphic to $Y^*(\frak{sl}_N)$. 
\endproclaim

\demo{Proof} 
It is easy to see that
$U_0(R_Y)$ is a graded quantization of $\frak{sl}_N(r)$, i.e.
it admits a $\Z$-grading (with $\deg(h)=1$) 
whose quasiclassical limit is the standard grading on $\g(r)$.
The Hopf algebra $Y^*(\frak{sl}_N)$ has the same property, as it is dual 
to the graded Hopf algebra $Y(\frak{sl}_N)$ with opposite coproduct. 
According to Drinfeld \cite{Dr1},
the graded quantization is unique. Therefore, Proposition 3.13 is proved. 
$\square$\enddemo

One can also consider the Yangian $Y(\frak{gl}_N)$. It is the 
quantized universal enveloping algebra with generators $t^*(u)=
t^*_{-1}u^{-1}+t^*_{-2}u^{-2}+...$, $t^*_{-i}\in Mat_N(k)\o Y(\frak{gl}_N)$, 
and defining relations (3.5),(3.6) (with $T^*$ instead of $T$), 
where $R=R_Y$, and $T^*(u)=
1+ht^*(u)=1+T^*_{-1}u^{-1}+
T^*_{-2}u^{-2}+...$. Let $Y^*({\frak {gl}}_N)$ be the dual Yangian, constructed
as in Section 3.1, with opposite product. 

\proclaim{Proposition 3.14} $U(R_Y)$ is isomorphic to $Y^*(\frak{gl}_N)$. 
\endproclaim

\demo{Proof} Let $H$ be the universal enveloping algebra of the 
abelian Lie algebra $t^{-1}k[t^{-1}]$. 
It is easy to see that we have Hopf algebra 
isomorphisms $U(R_Y)=U_0(R_Y)\o H$, $Y^*({\frak {gl}}_N)=
Y^*({\frak {sl}}_N)\o H$. 
Thus, Proposition 3.14 follows from Proposition 3.13.
$\square$\enddemo

\subhead 3.4. The PBW theorem for R-matrix algebras\endsubhead

In this section we prove
Proposition 3.9. 
Note that part (b) of Proposition 3.9 
plays the role of a Poincare-Birkhoff-Witt
theorem for R-matrix algebras. 

{\bf (a)} Let $Y=Y_r^{-1}Y_l$, where $Y_r,Y_l$ are the right and left 
hand sides of (3.7). Then from (3.5) we get
$$
[Y\o 1,T^{14}(u_1)T^{24}(u_2)T^{34}(u_3)]=0.
$$

Let $Y=\sum_{m\ge 0}Y_mh^m$. We prove that $Y_m$ is a scalar by induction.
The base of induction follows from the fact that $Y_0=1$.  
Assume that we proved that $Y_0,...,Y_{m-1}$ are scalars. 
If $F(R)$ is a flat deformation then, 
dividing the above equation by $h^m$ and reducing it mod h, we get
$[Y_m,g(u_1)\o g(u_2)]=0$ for any $g\in GL_N[[t]](k)$.
This implies that $Y_m$ is a scalar.   

Thus, we have shown that $Y$ is a scalar, which equals 1 mod h. 
Taking the determinants of both sides of
(3.7), we get $Y=1$. 

Similarly one shows that $R^{21}(-u)R(u)$ is a scalar.
Indeed, let $Z=R^{21}(v-u)R^{12}(u-v)$. Then from (3.5) we get
$$
[Z\o 1,T^{13}(u)T^{23}(v)]=0.
$$

Let $Z=\sum_{m\ge 0}Z_mh^m$. We prove that $Z_m$ is a scalar by induction.
The base of induction follows from the fact that $Z_0=1$.  
Assume that we proved that $Z_0,...,Z_{m-1}$ are scalars. 
If $F(R)$ is a flat deformation then, 
dividing the above equation by $h^m$ and reducing it mod h, we get
$[Z_m,g(u_1)\o g(u_2)\o g(u_3)]=0$ for any $g\in GL_N[[t]](k)$.
This implies that $Z_m$ is a scalar.   

Thus, we have shown that $Z$ is a scalar which equals 1 mod h, as desired.

{\bf (b)} First consider the case when 
$R=R_{rat}=1-\frac{h(\sigma-1/N)}{N(u-h/N^2)}$. 
In this case the flatness of $F(R)$ follows from the flatness 
of  $F(R_Y)$, since we have $R_{rat}=\frac{u}{u-h/N^2}R_Y$.

Now let $R$ be an arbitrary solution of (3.7), which satisfies
the assumption of (b). Let $R_s(u,h)=R(su,sh)$. 
Consider the algebra $F(R_s)$.

Assume the opposite, i.e. that $F(R)$ is not flat. 
Then $F(R_s)$ is not flat for any $s\ne 0$
(since $F(R_s)$ for $s\ne 0$ are all essentially the same). Denote by 
$F_\tau(R)$ 
the algebra $F(R_s)$ with $s$ being a formal 
parameter $\tau$. This is an algebra over $k[[h,\tau]]$, which is not a 
flat deformation of $k[GL_N[[t]]]$.  

Let $\Cal B$ be a basis of $k[GL_N[[t]]]$ consisting of monomials 
in the generators $t_{ij}^l$ (one has to choose an appropriate 
subset of the set of all monomials). It is clear that $\Cal B$ 
can be naturally regarded as a spanning system for $F_\tau(R)$, 
but it is linearly dependent. Pick a nontrivial linear relation 
$\sum_{j=1}^m \alpha_j(h,t)b_j=0$, where 
$b_j\in B$. 

Now we will make use of the representation theory of the dual Yangian. 
Let $V(a)$, $a\in k^*$, be the representation of $F(R_Y)$
on the vector space $V=k^N$  defined by 
$T(u)\to R_Y(u-a)$ (the shifted basic representation). 
We will consider all possible tensor products   
of these representations. 

\proclaim{Proposition 3.15} The natural map 
$F(R_Y)\to \oplus_n\oplus_{a_1,...,a_n\in k^*}\End(V(a_1)\o...\o V(a_n))[[h]]$
is injective. 
\endproclaim

\demo{Proof} It is enough to show that the map $U(R_Y)\to 
\oplus_n\oplus_{a_1,...,a_n\in k^*}\End(V(a_1)\o...\o V(a_n))[[h]]$
is injective. This follows easily from the fact that $U(R_Y)$ is
a deformation of $U(t^{-1}\g[t^{-1}])$. 
$\square$\enddemo

Now let $R$ be as in (b). Then, since 
$R$ satisfies (3.7), one can define the representations $V(a)$ of $F_\tau(R)$ 
in the same way as for the Yangian. Choose a finite collection of 
representations $W_1,...,W_l$ of the form $V(a_1)\o...\o V(a_n)$ 
such that $b_i$ are linearly independent in 
$M=(\End(W_1)\oplus...\oplus End(W_l))\o k[[h,\tau]]$.  
Let $P=\oplus_{j=1}^m kb_j\o_kk[[h,\tau]]$. 
We have a linear mapping $\theta: P\to M$, which assigns to $b_i$ the 
direct sum of the corresponding operators. We know that: 

(i) $P,M$ are finitely generated, free $k[[h,\tau]]$ modules. 
(ii) The map $\theta$ is not injective.
 (iii) $\theta$ is injective modulo $\tau$
(by Proposition 3.15).  

This is a contradiction,  
since by (ii) all the maximal minors of the matrix of $\theta$ in any basis 
are $0$, which contradicts (iii).

\subhead 3.5. $\d$-copseudotriangular structure on $F(R)$\endsubhead

In this section we describe a $\d$-copseudotriangular structure on 
any R-matrix algebra.

Let $\d$ be the derivation of $F(R)$ defined by the rule 
$\d T(u)=\frac{dT(u)}{du}$. 

\proclaim{Proposition 3.16} (i) Let $R$ satisfy (3.7). 
Then there exists a unique $\d$-copseudotriangular 
structure $B$ on $A(R)$ such that 
$$
B(T^{13}(u), T^{23}(v))(y)=R^{12}(u-v+y).\tag 3.13
$$

(ii) The form $B$ is given by the formula
$$
\gather
B(T^{1,p+q+1}(u_1)\dots T^{p,p+q+1}(u_p),
T^{p+1,p+q+1}(v_1)\dots T^{p+q,p+q+1}(v_q))(y)=\\
\prod_{i=1}^p\prod_{j=q}^1 R^{i,p+j}(u_i-v_j+y)
\tag 3.14
\endgather
$$
\endproclaim

\demo{Proof} Consider the form $B$ on the free algebra 
generated by $T(u)$ which is defined by (3.14).
It follows from the Yang-Baxter equation for $R$ that
the ideal generated by relations (3.5) is annihilated by $B$ on the 
right and on the left. This implies that $B$ descends to a form
on $F(R)$ and defines a form on $A(R)$.  
The fact that $B$ is a $\d$-copseudotriangular structure is easily obtained 
from the definition of $B$ and the properties of $R$. 
Uniqueness of $B$ is easily obtained
from equations (1.17). 
$\square$\enddemo

Let $f(u)\in 1+hk((u))[[h]]$, and 
$$
R^f(u)=f(u)R(u).\tag 3.15
$$ 
Then $R_Y^f(u)$
satisfies (3.7). 
Thus we obtain
the following proposition.

\proclaim{Proposition 3.17} For any $f\in 1+hk((u))[[h]]$,
there exists a unique $\d$-copseudotriangular
structure on $F(R)$ such that
$$
B_f(T^{13}(u), T^{23}(v))(y)=(R^f)^{12}(u-v-y).\tag 3.16
$$
\endproclaim

Let $f_0(u)$ be defined by the condition that $\qdet(R_Y^{f_0}(u))=1$. By
Proposition 3.17, 
$f_0$ defines a $\d$-copseudotriangular structure $B_{f_0}$ on
$Y(\frak{gl}_N)^\prime_{op}$. Since $\qdet(R_Y^{f_0}(u))=1$,
this $\d$-copseudotriangular structure descends to one on 
$Y(\frak{sl}_N)^\prime_{op}$. 
Thus, using Proposition 3.2, 
and the fact that $\d=d$ in this case, we get the following proposition. 

\proclaim{Proposition 3.18} 
There exists a unique copseudotriangular
structure on $Y(\frak{sl}_N)^\prime_{op}$ such that
$$
B(T^{13}(u), T^{23}(v))(y)=(R_Y^{f_0})^{12}(u-v+y).\tag 3.17
$$
This structure coincides with the one defined by (3.4). 
\endproclaim

\subhead 3.6. Elliptic R-matrix algebras\endsubhead

Here we consider R-matrix algebras associated to elliptic solutions
of the quantum Yang-Baxter equation.

Let $\Sigma$ be an elliptic curve over $k$, and $k(\Sigma)$ be the field
of rational functions on $\Sigma$. Let $du$ be an invariant 
differential on $\Sigma$, and $u$ the corresponding formal parameter 
near the origin. This parameter defines an embedding $k(\Sigma)\to k((u))$. 
Let $\Gamma$ be the group of points of order $N$ on $\Sigma$. 

\proclaim{Proposition 3.19} 
There exists a unique element $R(u,h)\in 
\End(k^N\o k^N)\o k(\Sigma)[[h]]$, 
such that $R=1+O(h)$, 
satisfying the following conditions.

(i) Let $R(u+h/N^2,h)=\sum_{m\ge 0}\rho_m(u)h^m$. Then $\rho_m(u)$ 
have at most simple poles at points of $\Gamma$, and no other 
singularities. 

(ii) Quasiperiodicity: if $\gamma\in\Gamma$, then
$R(u+\gamma)=(1\o\l(\gamma))(R(u))=
(\l(\gamma)^{-1}\o 1)(R(u))$, where $\l$ is as in Example 3 of 
Section 2.7.

(iii) $R(0):=\sum_{m\ge 0}\rho_m(-h/N^2)h^m=N\sigma$,
where $\sigma$ is the permutation matrix. 
\endproclaim

\demo{Proof} 

{\it Existence.} 
Let $r_\ell(u)$ be the elliptic R-matrix on $\Sigma$ introduced in 
Example 3 of Section 2.7. For $\g=\frak{sl}_N$, the
Casimir tensor $\Omega$ introduced in Chapter 2, is equal to
$\frac{\sigma-1/N}{N}$, so $r_\ell(u)=\frac{\sigma-1/N}{Nu}+O(1)$. 
Set 
$$
R(u,h)=N\sigma-h[r_\ell(u-h/N^2)-r_\ell(-h/N^2)].\tag 3.18
$$ 
It is easy to check that $R(u,h)$ satisfies conditions (i)-(iii). 

{\it Uniqueness.} Let $R_1,R_2$ be two functions satisfying
(i)-(iii). Let $R_1-R_2=h^mX+O(h^{m+1})$, $X\in\End(k^N\o k^N)\o k(\Sigma)$.

 From properties (i)-(iii) of $R_i$ we get
the following properties of $X$, respectively:

(i)' $X(u)$ has at most simple poles at $\Gamma$, and no other singularities. 

(ii)' $X(u+\gamma)=(1\o\l(\gamma))(X(u))=
(\l(\gamma)^{-1}\o 1)(X(u))$, $\gamma\in\Gamma$.

(iii)' $X(u)$ is regular at $u=0$. 

It is easy to see that any element $x\in \frak{gl}_N$ which is invariant under
the operators $\l(\gamma)$ is a scalar. Therefore, 
properties (i)'-(iii)'  
imply that $X$ is a scalar constant.  
Let $R_1-R_2=1+h^mX+h^{m+1}Z(u)+O(h^{m+2})$. 
The function $Z(u)$ satisfies properties (i)',(ii)'. Let $Z_0$ be the 
residue of $Z(u)$ at $u=0$. From condition (iii) we get $X-N^2Z_0=0$.
On the other hand, from conditions (i)',(ii)' for $Z(u)$ we get
$\sum_{\gamma\in\Gamma}(\l(\gamma)\o 1)(Z_0)=0$. This implies 
that $Z_0=X=0$, as desired. 
$\square$\enddemo

It is easy to check that the free term in the Laurent expansion of $r_\ell(u)$ 
vanishes. This implies that $R(u,h)=1-hr_\ell(u)+O(h^2)$. 
We will denote 
$R$ by $R_{\ell}$.

\proclaim{Proposition 3.22} $R_\ell(u,h)$ satisfies the quantum 
Yang-Baxter equation (3.7). 
\endproclaim

\demo{Proof} See \cite{Be1,Ch2,Ta}.
$\square$\enddemo

\proclaim{Corollary 3.23} The algebras $F(R_\ell)$, $A(R_\ell)$, $U(R_\ell)$ 
are flat deformations.
The algebra $U(R_\ell)$ is a quantization of the Lie bialgebra
$\frak{gl}_N(r_\ell)$.
\endproclaim

\demo{Proof} Follows from Proposition 3.9.
$\square$\enddemo

\subhead 3.7. Trigonometric R-matrix algebras\endsubhead

In this section we 
consider the limit of the R-matrix of Section 3.6 when the elliptic curve
degenerates into the multiplicative group $\Bbb G_m$.

Let $\g=\frak{sl}_N$. 
Let $\Sigma=\Bbb G_m$, $\Gamma$ be the group of points of order $N$, and
$\l: \Gamma\to \text{Aut}(\g)$ be the homomorphism defined by the rule
$\l(\e)=\Ad(\text{diag}(1,\e^r,...,\e^{r(N-1)}))$, where $1\le r\le N$ 
and $r$ is relatively prime to $N$.

Let $r_\tr(u)$ be the rational function on $\Sigma$ 
with values in $\g\o\g$ such that 
$r_\tr(u)=\frac{\Omega}{u}+O(1)$, $u\to 0$, 
$r_\tr$ is regular in $\bar\Sigma\setminus\Gamma$, where 
$\bar\Sigma=\Bbb P^1$ is the projective closure of $\Sigma$,
and $r_\tr(u+\gamma)=(1\o\l(\gamma))(r_\tr(u))$ for $\gamma\in\Gamma$.
It is clear that such a function is unique. 
It is obtained from $r_\ell(u)$ (for a suitable choice of $\l$) 
by the limiting procedure in which the elliptic curve degenerates 
into a rational curve. 

It is easy to write an explicit formula for $r_\tr$ using the global 
coordinate on $\Sigma$. This formula can be found in \cite{BD}.  
Namely, $r_\tr$ can be obtained from the r-matrix (2.13) by 
the transformation described in Proposition 2.5. 

Now we consider the quantum R-matrix. Define the function
$$
R_\tr(u,h):=N\sigma-h[r_\tr(u-\frac{h}{N^2})-r_\tr(-\frac{h}{N^2})].\tag 3.19
$$
As this function is a limiting case of $R_\ell(u,h)$,
Proposition 3.22 holds for $R_\tr$. Therefore, we have

\proclaim{Corollary 3.24} The algebras $F(R_\tr)$, 
$A(R_\tr)$, $U(R_\tr)$ are flat deformations.
The algebra $U(R_\tr)$ is a quantization of the 
Lie bialgebra
$\frak{gl}_N(r_\tr)$.
\endproclaim

\subhead 3.8. The factored Hopf algebras $F(R)_\z$, $U(R)_\z$\endsubhead

Here we quantize the Lie bialgebra $\g(r)_\z$, where 
$\g=\frak{gl}_N$, and $r$ is a classical r-matrix
which has a quantization $R$. 

 Let $\Sigma$ be a 1-dimensional algebraic group over $k$, and
$u$ be a canonical formal parameter near the origin.  
Let $R\in End(k^N\o k^N)(\Sigma)[[h]]$ be a function of the form 
$R=1-hr(u)+O(h^2)$ which satisfies the conditions of Proposition 3.9 (b). 
Let $A=A(R)$ be the h-formal group corresponding to $R$, then $U_A=U(R)$.
Let $\z=(z_1,...,z_n)$, $z_i\in\Sigma$, and $R(u)$ is regular 
at $z_i-z_j$ for $i\ne j$. In this case we can define 
the factored Hopf algebras $F(R)_\z,
A(R)_\z, U(R)_\z$ using the $\d$-copseudotriangular
structure on $A(R)$ defined in Section 3.5. 

\proclaim{Proposition 3.25} (i) The Hopf algebra $F(R)_\z$ 
is isomorphic to the h-adic completion of the 
algebra over $k[[h]]$ whose generators are the entries
of formal series $T_i(u)^{\pm 1}\in \End(k^N)\o F(R)_\z[[u]]$,
$T_i(u)=T^i_0+T^i_1u+...$,
and the defining relations are
$$
\gather
T_i(u)T_i(u)^{-1}=T_i(u)^{-1}T_i(u)=1,\\
R^{12}(u-v+z_i-z_j)T_i^{13}(u)T_j^{23}(v)=T_i^{23}(v)T_j^{13}(u)R^{12}(u-v
+z_i-z_j),
\tag 3.20\endgather
$$
and the coproduct is defined by
$$
\Delta(T_i(u))=T_i^{12}(u)T_i^{13}(u).\tag 3.21
$$

(ii) The algebra $F(R)_\z$ is a quantization of the Poisson group $G(r)_\z$. 
\endproclaim

\demo{Proof} By the definition, we have Hopf subalgebras 
$F(R)_{z_i}\subset F(R)_\z$ for all $i$, and multiplication
induces an isomorphism $F(R)_{z_1}\o...\o F(R)_{z_n}\to F(R)_\z$.
Let $T_i(u)$ is the image of the generating series $T(u)$ of
$F(R)_{z_i}$ under the  embedding $F(R)_{z_i}\to F(R)_\z$. 
Using Proposition 3.16 and the definition of $F(R)_\z$, it is easy 
to show that $T_i(u)$ satisfy relations (3.20).
Thus, we have a surjective homomorphism 
from the algebra defined by (3.20) to $F(R)_\z$. 
This homomorphism is an isomorphism modulo $h$, so it is an 
isomorphism.

Relation (3.21) is obvious. 
The fact that $F(R)_\z$ is a quantization of $G(r)_\z$ follows from the fact
that $F(R)_\z$ is a flat deformation. The Proposition is proved. 
$\square$\enddemo

\proclaim{Corollary 3.26} 
(i) The Hopf algebra $U(R)_\z$ 
is isomorphic to the h-adic completion of the 
algebra over $k[[h]]$ whose generators are the entries
of formal series $t_i(u)^{\pm 1}\in \End(k^N)\o U(R)_\z[[u]]$,
$t_i(u)=t^i_0+t^i_1u+...$,
and the defining relations are
$$
\gather
[t_i^{13}(u),t_j^{23}(v)]=
[r_*^{12}(u-v+z_i-z_j),t_i^{13}(u)+t_j^{23}(v)]+\\
h(r_*^{12}(u-v+z_i-z_j)t^{13}_i(u)t^{23}_j(v)-
t^{23}_j(v)t^{13}_i(u)r_*^{12}(u-v+z_i-z_j)),\tag 3.22\endgather
$$
where $r_*=h^{-1}(1-R)$, 
and the coproduct is defined by
$$
\Delta(t_i(u))=t_i^{12}(u)+t_i^{13}(u)+
ht_i^{12}(u)t_i^{13}(u).\tag 3.23
$$ 

(ii) The algebra $U(R)_\z$ is a quantization of the Lie bialgebra $\g(r)_\z$. 
\endproclaim

\head Appendix A: calculation of the square of the antipode\endhead

In this appendix we calculate the square of the antipode $S$ of a quantized
universal enveloping algebra $U_h(\g_+)$ obtained by quantization of
a Lie bialgebra $\g_+$ via the procedure of \cite{EK1}, Chapter 7-9.  
We will freely use the notation from \cite{EK1}.

Let $\g_+$ be a Lie bialgebra, and $\g$ the double of $\g_+$.
Let 
$M_-$, $M_+$ be the Verma modules for $\g$, and $M_+^*$ the dual 
Verma module (see Chapter 7 of \cite{EK1}). 

Let $1_-,1_+$ be the generating vectors in the Verma modules $M_+$,$M_-$ 
defined in Section 7.5 of \cite{EK1}, and $1_+^*\in M_+^*$ be the 
$\g$-invariant functional on $M_+$ normalized 
by $1_+^*(1_+)=1$. Let $(M_+^*)_1$ be the 
orthogonal complement to $1_+$ in $M_+^*$ (see the Appendix of \cite{EK1}). 

Let $\Cal M^e$ be the Drinfeld category for $\g$, defined in Chapter 7 
of \cite{EK1}.
Let $F$ be the fiber functor on $\Cal M^e$ defined in Chapter 
8 of \cite{EK1}, given by $F(V)=\Hom_{\Cal M^e}(M_-,M_+^*\o V)$.
Let $U_h(\g):=\End(F)$. 
For any $V\in\Cal M^e$, we identify $F(V)$ and $V[[h]]$ using the map
$\xi_V:F(V)\to V[[h]]$ defined by 
$\xi(a)=(1_+\o 1)(a1_-)$. Set $M_-^h=F(M_-)$, $M_+^{*h}=F(M_+^*)$, and 
denote by $1_+^{*h}\in F(M_+^*)$, $1_-^h\in M_-^*$,
$(M_+^{*h})_1\subset M_+^{*h}$ the images of $1_+^*,1_-,(M_+^*)_1$ 
under $\xi^{-1}$. Let $1_+^h: M_+^{*h}\to k[[h]]$ be the 
linear function which vanishes on $(M_+^{*h})_1$ and $1_+^h(1_+^{*h})=1$. 

Recall that the quantization $U_h(\g_+)$ of $\g_+$ defined in Chapter 9 of 
\cite{EK1} is equal, as a vector space, to
$F(M_-)=\Hom(M_-,M_+^*\o M_-)$.

\proclaim{Proposition A1} For any $a\in U_h(\g_+)$, 
$S^2(a)=\gamma^2\circ a$, where $\gamma$ is the braiding in $\Cal M^e$
defined in Section 7.7. of \cite{EK1}. 
\endproclaim

\demo{Proof} 
Let $R\in\End(F\times F)$ be the
universal $R$-matrix of $U_h(\g)$, i.e. $F(\gamma)=(\sigma R)^{-1}$, 
where $\sigma$ is
a permutation. We can represent $R$ as an infinite sum $R=\sum a_\alpha\o 
b_\alpha$, $a_\alpha\in U_h(\g_+),b_\alpha\in U_h(\g_-)$, where 
$\g_-$ is the dual algebra to $\g_+$. This series is convergent in the 
appropriate topology. Set $u:=\sum S^{-1}(a_\alpha)b_\alpha
\in U_h(\g)$ (it is easy to see that this series is also convergent). 

The following properties of $u$ were found by Drinfeld \cite{Dr3}:

(i) $u^{-1}xu=S^2(x)$ in $U_h(\g)$, where $x\in U_h(\g_+)$.

(ii) $\Delta(u)=(u\o u)(R^{21}R)^{-1}$. 

Now we are ready to prove the proposition. Since the functor $F$ is faithful,
it is enough to check the identity 
$$ 
F(S^2(a))=F(\gamma^2)\circ F(a) \text{ in }\Hom_{U_h(\g)}(M_-^h,
M_+^{*h}\o M_-^h).\tag A1
$$ 
Applying the transformation $x\to (1_+^{h}\o 1)(x1_-^h)$ to 
both sides of (A1), and using the fact that $F(\gamma^2)=(R^{21}R)^{-1}$,
we obtain an equivalent identity
$$
S^2(a)=(1_+^{h}\o 1)((R^{21}R)^{-1} F(a)1_-^h)\text{ in }M_-^h.\tag A2
$$
 Using property (ii) 
of $u$, we can rewrite this identity 
in the form
$$
S^2(a)=(1_+^{h}\o 1)((u^{-1}\o u^{-1})\Delta(u) F(a)1_-^h).\tag A3
$$
Since $(\e\o 1)(R)=(1\o \e)(R)=1$,
we see that $u1_-^h=1_-^h$, $u(M_+^*)_1=(M_+^*)_1$. 
Using these properties of $u$, we reduce (A3) to the form
$$
S^2(a)=(1_+^{h}\o 1)((1\o u^{-1})F(a)1_-^h).\tag A4
$$
The right hand side of (A4) equals $u^{-1}*a$, where $u^{-1}*v$ denotes the 
action of
$u^{-1}$ on $v\in F(V)$. Thus, we have reduced (A2) to
$$
S^2(a)=u*a.\tag A5 
$$
 
On the other hand, from property (i) of $u$, and the equation $u1_-^h=1_-^h$,
we get $S^2(a)1_-^h=uau^{-1}1_-^h=ua1_-^h=u*a$, which coincides with (A5). 
The proposition is proved. 
$\square$\enddemo

Let $g_\alpha$ be a basis of $\g_+$, $g_\alpha^*$ the dual basis of $\g_-=
\g_+^*$, and $C_+=\sum_\alpha g_\alpha g_\alpha^*$ be the half-Casimir. 
$C_+$ is an endomorphism of the forgetful functor 
to $k[[h]]$-modules on the category $\Cal M^e$. 

\proclaim{Proposition A2} We have 
$$
\xi_{M_-}\circ S^2 =
e^{hC_+}\circ \xi_{M_-}.\tag A6
$$
\endproclaim

\demo{Proof} We start with 
a tautological identity
$$
(1_+\o 1)((1\o e^{hC_+})a1_-)=e^{hC_+}(1_+\o 1)(a1_-).\tag A7
$$
Let $\Omega$ be the Casimir operator for $\g$ (see Chapter 7 of \cite{EK1}).
Since $\Omega=\Delta(C_+)-C_+\o 1-1\o C_+$, and $C_+1_-=0$, $S(C_+)1_+=0$,
formula (A7) implies
$$
(1_+\o 1)(e^{-h\Omega}a1_-)=e^{hC_+}(1_+\o 1)(a1_-).\tag A8
$$
By Proposition A1, $S^2(a)=\gamma^2\circ a$. 
By the definition, in $\Cal M^e$ we have $\gamma^2=e^{-h\Omega}$, 
Thus, (A8) reduces to the form
$$
(1_+\o 1)(S^2(a)1_-)=e^{hC_+}(1_+\o 1)(a1_-),\tag A9
$$
which is equivalent to (A6).
The proposition is proved.
$\square$\enddemo

\proclaim{Proposition A3} Let $d$ be the derivation of $\g_+$ given 
by the formula $d=-\frac{1}{2}\mu\circ \delta$, where $\mu,\delta$ are 
the commutator and cocommutator in $\g_+$. Let
$D=\frac{1}{h}\ln S^2$ be the canonical derivarion of $U_h(\g_+)$. 
Then the image of $d$ under the functor $U_h$ equals to $D$. 
\endproclaim

\demo{Proof} We first prove the following Lemma.

{\bf Lemma.} For any $x\in \g_+$, $dx=[C_+,x]$ in any $V\in \Cal M^e$.    

{\it Proof of the Lemma.} 
$$
[C_+,x]=\sum (g_\alpha[g_\alpha^*,x]+[g_\alpha,x]g_\alpha^*)\tag A10
$$
For any $y^*\in\g_-$, we have $[y^*,x]=(1\o y^*)(\delta(x))-ad^*x(y^*)$. 
Using this identity and the invariance of the pairing
between $\g_+,\g_-$ under $x$, we get
$$
[C_+,x]=\sum g_\alpha(1\o g_\alpha^*)(\delta(x))=m_{21}(\delta(x)),\tag A10
$$
where $m_{21}(a\o b):=ba$. Since $\delta(x)$ is skew symmetric, 
we have $m_{21}(\delta(x))=-\frac{1}{2}\mu(\delta(x))=dx$. The Lemma is 
proved. 

Now we prove Proposition A3. It follows from proposition A2 that after 
$U_h(\g_+)$ is identified with $M_-[[h]]=U(\g_+)[[h]]$, as a vector space,
via $\xi_{M_-}$, the derivation $D$ is given by $Dx=[C_+,x]$. 
Therefore, by the lemma, $D=d$, as desired.
$\square$\enddemo

\end
\Refs

\ref\by [AGS1] Alekseev, A., Grosse, H., and Schomerus, V. 
\paper Combinatorial quantization of the Hamiltonian Chern-Simons theory I
\jour CMP\vol 172\yr 1995\endref

\ref\by [AGS1] Alekseev, A., Grosse, H., and Schomerus, V. 
\paper Combinatorial quantization of the Hamiltonian Chern-Simons theory II
\jour CMP\vol 174\yr 1995\endref

\ref\by [Be1] A.A.Belavin\paper Discrete groups and integrability 
of quantum systems\jour Funct. Anal. Appl.\vol 14\issue 4\yr 1980\endref

\ref\by [Be2] A.A.Belavin\paper Dynamical symmetry of integrable
quantum systems\jour Nucl. Phys.\vol B180(FS2)\issue 2\pages 189-200\yr
1981\endref

\ref\by [BD] A.A.Belavin and V.G.Drinfeld\paper Solutions of the classical
Yang-Baxter equations for simple Lie algebras\jour Funct. Anal. Appl.
\vol 16\pages 159-180\yr 1982\endref

\ref\by [Ch1] I.V.Cherednik\paper On R-matrix quantization of formal 
loop groups\jour Group Theoretical Methods in Physics, VNU Publ., 
\yr 1986\endref

\ref\by [Ch2] I.V.Cherednik\paper On properties of factorized S-matrices  
in terms of elliptic functions\jour Yad. Fiz. \vol 36\issue 2\pages 549-557
\yr 1982\endref

\ref\by [Dr1] V.G.Drinfeld \paper Quantum groups \jour 
Proc. Int. Congr. Math. (Berkeley, 1986)\vol 1\pages 798-820\endref 

\ref\by [Dr2] V.G.Drinfeld \paper Hopf algebras and the quantum 
Yang-Baxter equation \jour Soviet Math. Doklady\vol 32\yr 1985\endref

\ref\by [Dr3] V.G.Drinfeld \paper 
On almost cocommutative Hopf algebras \jour Len. Math.J.
\vol 1\pages 321-342\yr 1990\endref

\ref\by [EK1] P.Etingof and D. Kazhdan\paper Quantization of Lie bialgebras, I,
q-alg 9506005\jour Selecta math. \vol 2\issue 1\yr 1996\pages 1-41\endref

\ref\by [EK2] P.Etingof and D. Kazhdan\paper Quantization of Lie bialgebras, 
II, \yr 1996\jour q-alg 9701038\endref

\ref\by [FRT] N.Yu.Reshetikhin, L.A.Takhtajan, and L.D.Faddeev,  
\paper Quantization of Lie groups and Lie algebras\jour Len. Math. J.\vol. 
1\yr 1990\pages 193-225\endref

\ref\by [RS] N.Yu.Reshetikhin and M.A.Semenov-Tian-Shansky\paper
Central extensions of quantum current groups\jour Lett. Math. Phys.
\vol 19\yr 1990\pages 133-142\endref

\ref\by [Ta] L.Takhtajan \paper Solutions of the triangle 
equations with $\Z_n\times \Z_n$-symmetry and matrix analogs of the Weierstrass
zeta and sigma-functions\jour Zapiski Nauch. Sem. LOMI\vol 133\pages 
258-276\yr 1984\endref
